\documentclass[preprint,12pt]{elsarticle}

\usepackage{amssymb}
\usepackage{pifont}
\usepackage{amsmath}
\allowdisplaybreaks
\usepackage{amsthm}
\usepackage{graphicx}
\usepackage[normalem]{ulem}
\usepackage{url}
\usepackage{hyperref}

\usepackage{booktabs}  
\usepackage{longtable} 
\usepackage{tabularx}

\usepackage[acronym, toc, nonumberlist]{glossaries} 

\makeglossaries 

\newacronym{pwl}{PWL}{piecewise linear}
\newacronym{lp}{LP}{linear programming}
\newacronym{slp}{SLP}{sequential linear programming}
\newacronym{qp}{QP}{quadratic programming}
\newacronym{milp}{MILP}{mixed-integer linear programming}
\newacronym{nlp}{NLP}{nonlinear programming}
\newacronym{minlp}{MINLP}{mixed-integer nonlinear programming}
\newacronym{doe}{DOE}{Department of Energy}
\newacronym{nlr}{NLR}{National Laboratory of the Rockies}
\newacronym{pnnl}{PNNL}{Pacific Northwest National Laboratory}
\newacronym{huc}{HUC}{Hydro Unit Commitment}
\newacronym{poz}{POZ}{prohibited operating zones}
\newacronym{psh}{PSH}{pumped-storage hydropower}
\newacronym{relu}{ReLU}{Rectified Linear Unit}
\newacronym{sos2}{SOS2}{Special Ordered Set of type 2}
\newacronym{1d}{1D}{one-dimensional}
\newacronym{2d}{2D}{two-dimensional}

\journal{Applied Energy}

\begin{document}

\begin{frontmatter}



\title{Modeling Gaps in Hydropower Cascading System Models: A Systematic Review of Rule-Based Formulations}

\author[1]{Quentin Ploussard\corref{cor1}}
\ead{qploussard@anl.gov}

\author[1]{Lukas Livengood}

\author[2]{Slaven Kincic}

\cortext[cor1]{Corresponding author}

\affiliation[1]{organization={Argonne National Laboratory},
            addressline={9700 South Cass Avenue}, 
            city={Lemont},
            postcode={60439}, 
            state={IL},
            country={USA}}
\affiliation[2]{organization={Pacific Northwest National Laboratory},
            addressline={902 Battelle Boulevard}, 
            city={Richland},
            postcode={99352}, 
            state={WA},
            country={USA}}

\begin{abstract}
The coordination of cascading hydropower systems represents a fundamental challenge in modern energy systems engineering, requiring a sophisticated balance between multi-reservoir physics, stringent environmental regulations, and dynamic market participation. As intermittent energy sources increase, the transition to high-fidelity hydropower modeling has become a core requirement for ensuring power system reliability, long-term energy resilience and affordability. This review provides a comprehensive analysis of 131 seminal articles through the exclusive lens of optimization-based approaches, intentionally omitting pure simulation and heuristic methods to focus on rigorous mathematical formulations. A generalized 10-equation mathematical framework is established as a formal baseline, capturing the full physical and hydraulic behavior of cascading systems, including spatiotemporal inflow routing, storage-to-elevation relationships, head-dependent power generation, and prohibited operating zones. Each article is evaluated against every equation of this Standard Model in a systematic census documenting which physical relationships are included, simplified, or omitted, providing an empirical measure of modeling fidelity across the field. 
Modeling simplifications are evaluated through the lens of grid reliability rather than water management performance alone. The review makes a focused technical case for mixed-integer linear programming with piecewise linear approximations as the optimal balance between physical accuracy and computational tractability, highlighting recent optimal regression techniques that minimize combinatorial overhead. Finally, a bibliometric analysis of solver usage identifies the near-absence of high-performance open-source solvers as a critical reproducibility barrier, and a promising avenue for broader adoption of high-fidelity cascading hydropower models.
\end{abstract}


\begin{highlights}
\item \textbf{Systematic census:} Equation-level analysis of 131 cascading hydropower articles.

\item \textbf{Modeling gaps:} Nominal head assumptions and missing POZ obscure grid stability risks.

\item \textbf{Fidelity vs. tractability:} MILP with PWL approximations identified as the optimal strategy.

\item \textbf{Combinatorial efficiency:} Optimal PWL regression minimizes binary variables in MILP models.

\item \textbf{Open-source solvers:} HiGHS and SCIP absent from the literature despite competitive performance.
\end{highlights}


\end{frontmatter}




\printglossary[type=\acronymtype, title={Abbreviations}]

\section*{Nomenclature}

\addcontentsline{toc}{section}{Nomenclature} 

\subsection*{Sets and domains}
\begin{tabular}{lp{12cm}}
$\mathcal{R}$ & Set of reservoirs in the cascading system, indexed by $r$.\\
$\mathcal{P}$ & Set of hydropower plants or units, indexed by $p$.\\
$\mathcal{P}(r)$ & Set of hydropower plants or units of reservoir $r$.\\
$\mathcal{T}$ & Set of time periods in the optimization horizon, indexed by $t$.\\
$\mathcal{Z}_p$ & Set of operating zones for the power output $P_{p,t}$ and net hydraulic head $h_{p,t}$ of unit $p$.\\
$\mathcal{A}$ & Set of hydraulic arcs $r' \to r$ representing the topology of the cascading reservoir system.\\
$r(p)$ & Mapping of plant $p$ to its associated reservoir $r$.\\
$d(r)$ & Immediate downstream reservoir or afterbay of reservoir $r$.
\end{tabular}

\subsection*{Parameters}
\begin{tabular}{lp{12cm}}
$W_{r,t}$ & Natural local inflow to reservoir $r$ during period $t$ ($m^3/s$).\\
$\Delta t$ & Duration of each time step $t \in \mathcal{T}$ ($s$).\\
$\rho$ & Density of water ($kg/m^3$).\\
$g$ & Gravitational acceleration ($m/s^2$).\\
$\overline{E}_r, \underline{E}_r$ & Maximum and minimum allowable water surface elevation at reservoir $r$ ($m$).\\
\end{tabular}

\subsection*{Variables}
\begin{tabular}{lp{12cm}}
$\mathcal{F}$ & Generalized objective function.\\
$P_{p,t}$ & Power output of plant/unit $p$ during period $t$ ($MW$).\\
$V_{r,t}$ & Storage volume of reservoir $r$ at the end of period $t$ ($m^3$).\\
$E_{r,t}$ & Water surface elevation (forebay) of reservoir $r$ at the end of period $t$ ($m$).\\
$h_{p,t}$ & Net hydraulic head of unit $p$ during period $t$ ($m$).\\
$\eta_{p,t}$ & Efficiency coefficient of plant/unit $p$ during period $t$ ($p.u.$).\\
$Q_{r,t}$ & Total outflow from reservoir $r$ during period $t$ ($m^3/s$).\\
$Q^{\text{P}}_{p,t}$ & Turbined (power-generating) discharge through unit $p$ at period $t$ ($m^3/s$).\\
$Q^{\text{NP}}_{r,t}$ & Non-power discharge (spillage/bypass) from reservoir $r$ at period $t$ ($m^3/s$).\\
$I_{r,t}$ & Total inflow to reservoir $r$ at period $t$, including upstream releases ($m^3/s$).\\
$L_{r,t}$ & Water losses (evaporation/seepage) from reservoir $r$ at period $t$ ($m^3/s$).
\end{tabular}

\subsection*{Functions}
\begin{tabular}{lp{12cm}}
$F_r(\cdot)$ & Routing function of reservoir $r$.\\
$f_r^{\text{E}}(\cdot)$ & Storage-to-elevation functions of reservoir $r$.\\
$f_p^{\text{P}}(\cdot)$ & Power output function of unit $p$.\\
$f_r^{\text{TR}}(\cdot)$ & Tailrace effect function of reservoir $r$.\\
$f_p^{\text{HL}}(\cdot)$ & Head loss effect function of unit $p$.
\end{tabular}

\section{Introduction}

The development of high-fidelity modeling frameworks for cascading hydropower systems has become a critical priority for ensuring national energy resilience and affordability \cite{doe_reliability_topic}. This imperative is underscored by the U.S. \gls{doe} HydroWIRES initiative, which emphasizes the technical need to advance the representation of hydropower resources in electric power system models and to address the ``seams'' between water and power modeling domains \cite{osti_1726280}. As noted in foundational technical workshops conducted by \gls{pnnl} and the \gls{nlr}, bridging these gaps is essential to accurately characterize the flexibility and ``headroom'' of hydroelectric assets within large-scale interconnection studies \cite{osti_1726280, 11225846, 11225344}. A foundational technical assessment by Kincic et al. \cite{osti_1922507} identified and prioritized the most critical hydropower modeling gaps in power system planning and operational studies through direct engagement with Western Interconnection industry experts. Their findings confirm that hydro-based generation in current planning software accounts for nominal water availability only, and that interdependencies among cascading projects (including headwater elevation constraints, environmental flow requirements, and tailwater effects) are systematically absent from base cases, leading to unrealistic dispatch and overly optimistic reliability assessments.

Comprehensive power system models (commonly referred to as basecases) serve as the foundational tools for both operational studies (real-time to seasonal horizons) and long-term planning studies (capital investment and system expansion). Despite their distinct time horizons, both study types share the same steady-state and dynamic grid representations. While hydropower accounts for approximately 25\% of total generation in the Western Interconnection, current basecases concentrate exclusively on the electrical aspects of the system, systematically overlooking water availability, environmental flow requirements, and the interdependencies among cascading hydroelectric projects \cite{osti_1922507}. As the penetration of variable energy resources increases, this omission becomes increasingly consequential: unrealistic hydropower dispatch in basecase models can overestimate available flexibility and, in extreme cases, contribute to partial or complete blackouts. Accurately representing water management constraints within power system models is therefore essential for both reliable system operation and sound long-term planning, a gap that the present review directly addresses through the rigorous mathematical formalization of cascading hydropower interdependencies.

A primary motivation for utilizing rule-based mathematical programming over alternative paradigms is the requirement for a rigorous representation of physical and coupling constraints. While statistical and machine learning models have gained traction due to their computational speed, they possess significant disadvantages in operational planning. These models lack the explicit structural constraints defined by physical laws; consequently, they struggle to remain accurate if underlying operating rules change, \textit{e.g.}, new environmental flow regulations or altered market structures, as they often lack sufficient historical data to represent these novel regimes \cite{osti_1726280, 11225846}. 

Furthermore, the mathematical representation of a cascading system is significantly more complex than that of isolated hydropower plants. In isolated systems, the primary constraint is typically local reservoir volume; however, cascading systems are governed by strict interdependency rules where the discharge and spillage decisions of an upstream facility fundamentally determine the available potential of downstream units \cite{8633891, 7790824, TeegavarapuSimonovic2000}. The inclusion of these spatial and temporal interdependencies, such as inflow-discharge coupling, is vital; failing to model these ``ripple effects'' leads to inaccurate flexibility assessments that can threaten grid stability \cite{11225344, 11225846}.


The existing literature is supported by several comprehensive reviews that have approached reservoir and hydropower optimization through complementary but distinct lenses. Labadie \cite{labadie2004optimal} provided a foundational assessment of multireservoir system optimization, scrutinizing \gls{lp}, \gls{nlp}, dynamic programming, and early genetic algorithm methods, and identified a persistent gap between theoretical developments and real-world implementations driven by model complexity and operator skepticism. More recently, Lai et al. \cite{lai2022review} tracked the evolution of reservoir operation optimization from traditional \gls{lp} and dynamic programming models to modern metaheuristic algorithms (including evolutionary, swarm intelligence, and nature-inspired techniques) over the decade from 2011 to 2021, evaluating performance through reservoir system policy metrics such as reliability, resilience, and vulnerability. Complementarily, Parvez et al. \cite{w11071392} surveyed a broad range of hydro generation scheduling techniques across multiple time horizons, classifying methods into heuristic, mathematical programming, and hybrid categories, and concluded that \gls{milp} is most effective for large-scale systems while Lagrangian relaxation remains preferred for computational speed. Taktak and D'Ambrosio \cite{taktak2017overview} offered a more specialized review of mathematical programming approaches, covering \gls{milp}, \gls{minlp}, dynamic programming, Lagrangian relaxation, and Benders decomposition (specifically for the deterministic hydro unit commitment problem in hydro valleys), with a particular focus on linearization strategies for the nonlinear head-dependency of power output. From a broader system integration perspective, Thirunavukkarasu et al. \cite{THIRUNAVUKKARASU2023113192} evaluated optimization techniques for hybrid energy systems spanning intermittent energy sources, hydropower, and battery storage, noting that hydropower sources are consistently simplified or under-represented relative to variable energy resources, and that AI-based hybrid algorithms outperform classical methods for global optimization across these multi-component systems.

Despite the breadth of these reviews, a persistent ``modeling gap'' remains regarding the integration of cascading hydropower into large-scale grid reliability assessments. Existing works either treat hydropower plants as isolated or aggregated units, prioritize water management objectives over grid stability, rely on metaheuristic methods that cannot guarantee physical feasibility, or restrict their scope to single-valley unit commitment without addressing inter-reservoir coupling. The present article directly addresses these gaps through the following distinct contributions:
\begin{itemize}
    \item Standardized cascading model and taxonomic framework: A generalized 10-equation mathematical model is established as a formal baseline, capturing the full physical and hydraulic behavior of cascading systems, including spatiotemporal inflow routing, storage-to-elevation relationships, head-dependent power generation, and prohibited operating zones. This ``Standard Model'' serves as a structured reference for classifying the \gls{lp}, \gls{milp}, and \gls{minlp} adaptations employed across the 131 reviewed articles, providing a level of modeling detail absent from prior reviews.

    \item Equation-level systematic census of 131 articles: Each reviewed article is evaluated against every fundamental constraint of the Standard Model, documenting which physical relationships are included, simplified, or omitted, and how nonlinear functions are approximated. This per-article, per-equation census provides an empirical foundation for quantifying the current state of modeling fidelity across the field.

    \item Fidelity gap framing for grid reliability: The modeling simplifications identified in the literature, \textit{e.g.}, nominal hydraulic head assumptions and the omission of prohibited operating zones, are evaluated through the lens of power grid reliability rather than water management performance alone. This framing, aligned with the \gls{doe} HydroWIRES initiative, exposes stability risks that are systematically obscured by prevailing modeling conventions.

    \item \gls{milp} as a strategic compromise via \gls{pwl} universal approximation: A focused technical case is made for \gls{milp} with \gls{pwl} approximations as the optimal balance between physical fidelity and computational tractability. This analysis connects \gls{pwl} functions to their theoretical foundation as universal approximators, sharing the mathematical structure of \gls{relu}-based neural networks, and evaluates optimal regression techniques for minimizing combinatorial overhead in high-fidelity cascading models.

    \item Open-source solver gap and reproducibility roadmap: A bibliometric analysis of solver usage reveals that high-performance open-source solvers such as SCIP and HiGHS are virtually absent from the existing corpus despite demonstrated competitive performance. This finding identifies a critical barrier to broader adoption and reproducibility, and positions emerging open-source tools as a viable path toward accessible, high-fidelity cascading hydropower optimization.
\end{itemize}

The remainder of this article is organized as follows: 
Section~\ref{sec:scope} defines the scope of this systematic review and provides the detailed search methodology used to curate the 131-article corpus. 
Section~\ref{sec:bibliometrics} provides a bibliometric overview, highlighting the interdisciplinary breadth of journals, temporal publication trends, and the diversity of optimization objectives, modeling paradigms, and numerical solvers identified in the literature. Section~\ref{sec:math} establishes a generalized hydropower cascading model and evaluates the physical simplifications required to adapt this ``Standard Model'' into various mathematical formulations, including \gls{lp}, \gls{milp}, and \gls{minlp}. This section includes a comprehensive analysis of the specific model versions and coupling equations utilized across the reviewed studies. 
Section~\ref{sec:discussion} provides a critical discussion on the trade-offs between modeling fidelity and computational tractability. This section explores the strategic advantages of \gls{milp} formulations and \gls{pwl} approximations, specifically regarding strong mathematical formulations and optimal regression techniques, alongside the role of modern open-source solvers as a promising trajectory for high-fidelity cascading models. 
Section~\ref{sec:conclusion} summarizes the findings and offers concluding remarks on future research directions.

\section{Scope and methodology}

\label{sec:scope}

To evaluate the current state of cascading hydropower optimization, this review employs a systematic approach to literature selection and classification. This section defines the search strategy and the rigorous criteria used to curate the 131 analyzed works.

\subsection{Review methodology}
The literature search was conducted across established academic databases, including Scopus, IEEE Xplore, and Google Scholar. To identify emerging research trends and ensure comprehensive coverage, AI-powered discovery tools, specifically Consensus and Elicit, were also used.

The search query was constructed using technical keywords including: \textit{Cascade, Hydropower, Scheduling, Optimization, Metaheuristic, Linear, Nonlinear, Programming, Mixed Integer, LP, NLP, MILP, and MINLP}.

\subsection{Inclusion and exclusion criteria}
The selection process was specifically designed to isolate and evaluate the mathematical representations of hydropower cascading systems. To maintain this technical focus, the following exclusion criteria were applied:

\begin{itemize}
\item \textbf{Mathematical transparency:} To serve as a rigorous technical benchmark, articles that did not provide explicit physical equations, such as hydropower equations or cascading mass-balance equations, were excluded from the analysis.
\item \textbf{Methodological focus:} To ensure technical rigor, the literature
selection prioritized multi-reservoir cascading architectures over isolated,
single-reservoir models. The inclusion boundary was drawn at the level of
formulation transparency rather than solving method: an article was
retained if and only if it explicitly defines the physical and hydraulic
constraints of the cascading system in mathematical form, enabling a direct
comparison of modeling choices across the corpus. As a result,
metaheuristic-solved models, such as those employing Genetic Algorithms
or Particle Swarm Optimization, are retained when they provide explicit
cascading formulations, while pure search-based approaches that optimize
over simulation outputs without a defined constraint structure are excluded.
This boundary reflects the core objective of this review: to evaluate the
fidelity of physical representations, not to assess the relative performance
of optimization algorithms. Furthermore, as this study specifically aims to
bridge the operational ``seams'' between power systems and water management,
the selection prioritized cascading models involving hydropower-generating
facilities, with the analysis of each article focusing on the
power-generating components and their associated hydraulic coupling
constraints.
\item \textbf{Technology scope:} Although the standardized modeling framework introduced in this review does not explicitly incorporate the specialized operational constraints of \gls{psh} plants, literature involving systems with \gls{psh} assets was not excluded. In these instances, the analysis focuses strictly on the conventional cascading formulations and inter-reservoir dynamics rather than the specific pumping-cycle mechanics.
\end{itemize}

Ultimately, this filtering process resulted in a core dataset of 131 articles that directly address the ``modeling gaps'' in high-fidelity cascading hydropower dispatch.

\section{Bibliometric overview}
\label{sec:bibliometrics}

This section provides an overview of the 131 reviewed articles, detailing the bibliographic data and technical characteristics of the corpus. It maps the interdisciplinary distribution of publication venues and traces the temporal evolution of the field. Additionally, the overview categorizes the literature based on optimization objectives, mathematical modeling paradigms, and the landscape of commercial and open-source numerical solvers.

\subsection{Interdisciplinary publication landscape}
\begin{figure}[ht!]
\centering
\includegraphics[width=1.0\textwidth]{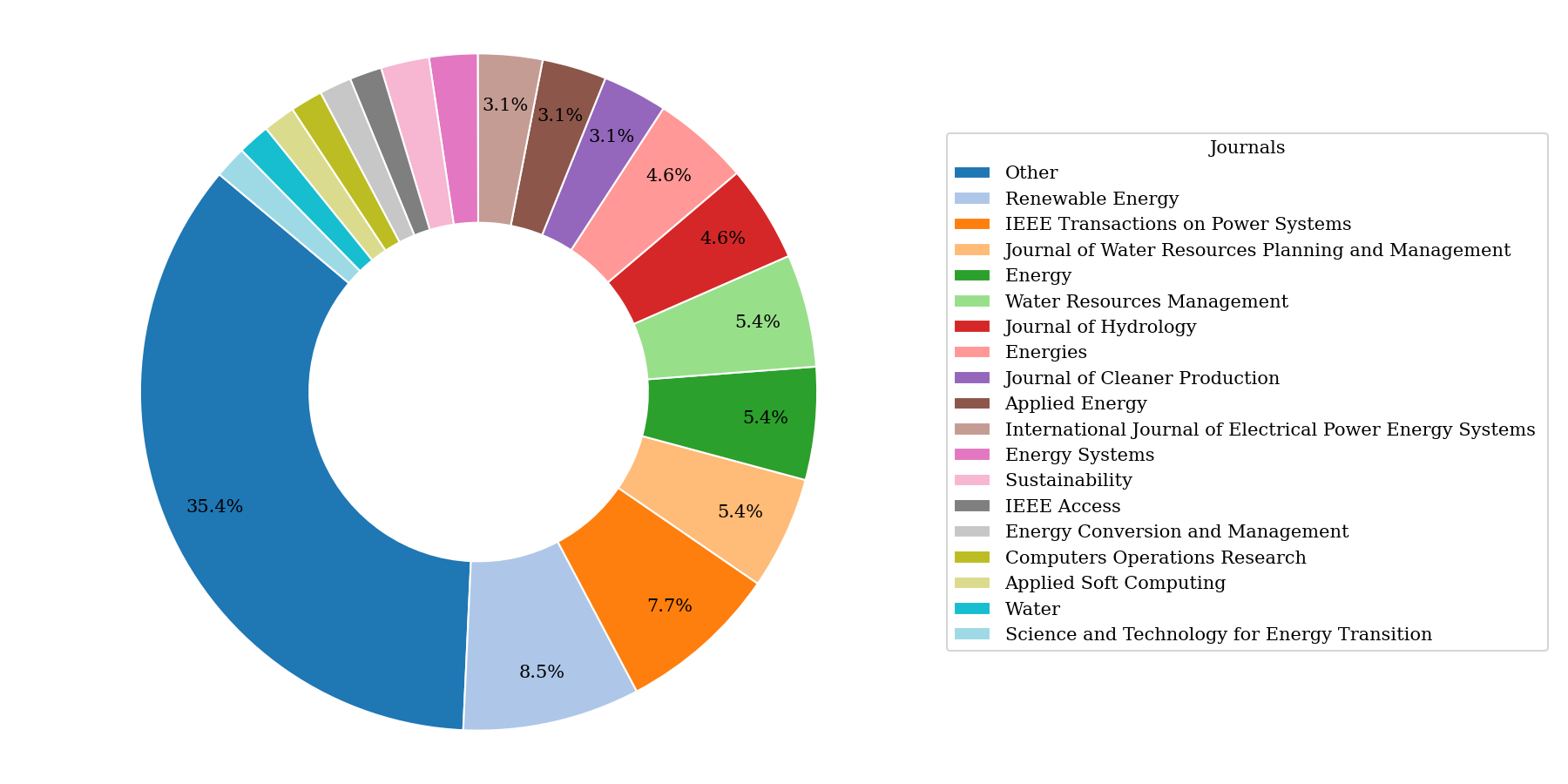} 
\caption{Publication distribution by venue}
\label{fig:pie_chart}
\end{figure}

The distribution of the 131 selected articles reveals a highly interdisciplinary research landscape, as illustrated in Fig.~\ref{fig:pie_chart}. The research topic is a significant area of interest for journals across water management \cite{TeegavarapuSimonovic2000, LU2021114055, SLP01, AmaniHUC, MMLPNiu, DoganMix}, hydrology \cite{SU2020125556, YOO2009182, CHEN2023129185, Wang2017OptOutputError, Zhu2025PFMODO, LU2021126388, HUANG2025132756}, power systems \cite{5565530, 9211795, 7790824, 1137622, 7038354, LIAO2024121341, SwarmMINLP}, and operations research \cite{7038354, SLP01, FLETEN20082656, NIU2021107315, AMINABADI2024282, NIU2018562}. This reflects the multiple interconnected roles served by water reservoirs in the United States, most notably water management and power generation, and highlights the importance of modeling these functions jointly rather than in isolation.

\subsection{Temporal evolution and objective trends}
\begin{figure}[ht!]
\centering
\includegraphics[width=1.0\textwidth]{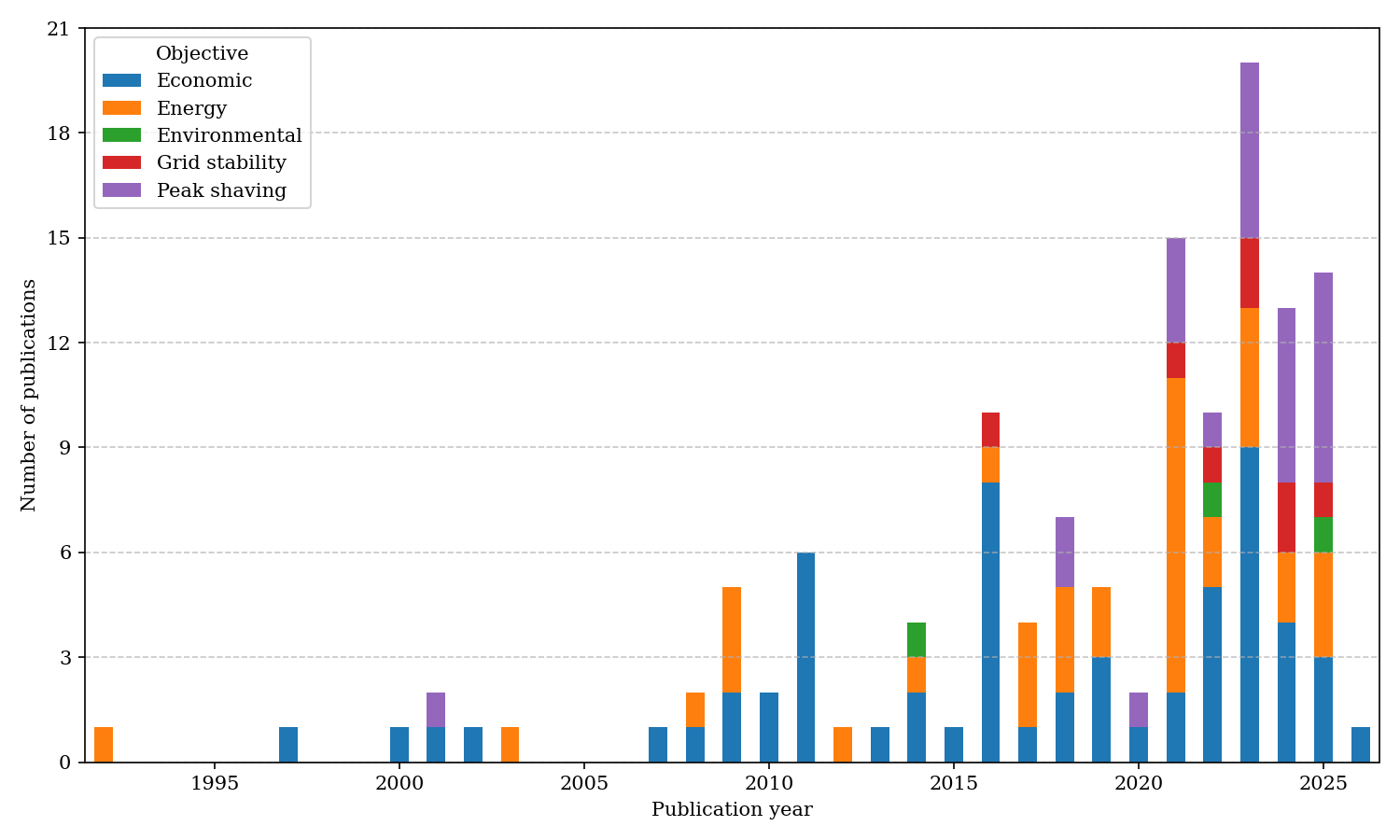} 
\caption{Yearly publication trends}
\label{fig:yearly_dist}
\end{figure}

The temporal distribution of the literature, shown in Fig.~\ref{fig:yearly_dist}, indicates a clear and accelerating interest in the modeling of cascading systems. This trend may reflect a growing focus on accurately representing these inherently complex systems, supported by advances in computational hardware and solver algorithms that make such modeling increasingly tractable.

A detailed analysis of the optimization objectives identified across the corpus reveals a growing focus on hydropower scheduling for peak shaving \cite{9211795, FANG2024120932, 8494460, SU2022395, en14040887, ZHAO2024129834} and grid stability purposes \cite{LU2024122085, WANG202268, LIAO2021970, LIAO2023127685, Liu2016_JointOptimization, WANG2024130258}, although economic \cite{8633891, 5565530, 1137622, TeegavarapuSimonovic2000, HOSNAR2014194, LIAO2024121341} and energy objectives \cite{7790824, 7038354, 9314115, 7447809, WU2023507, YOO2009182} remain an important goal of these models. This evolution is likely driven by the increasing penetration of variable energy sources, which requires more flexible, high-fidelity dispatch solutions to maintain system reliability. 

In contrast, relatively few studies explicitly focus on environmental objectives \cite{Forknall2014, Zhu2025PFMODO, CHENG2022123908}. This may be because environmental considerations are more often incorporated as operational constraints (\textit{e.g.}, target release volumes, flow rate limits, ramp limits) rather than being formulated as primary optimization objectives.

\subsection{Mathematical paradigms and formulations}
\begin{figure}[ht!]
\centering
\includegraphics[width=1.0\textwidth]{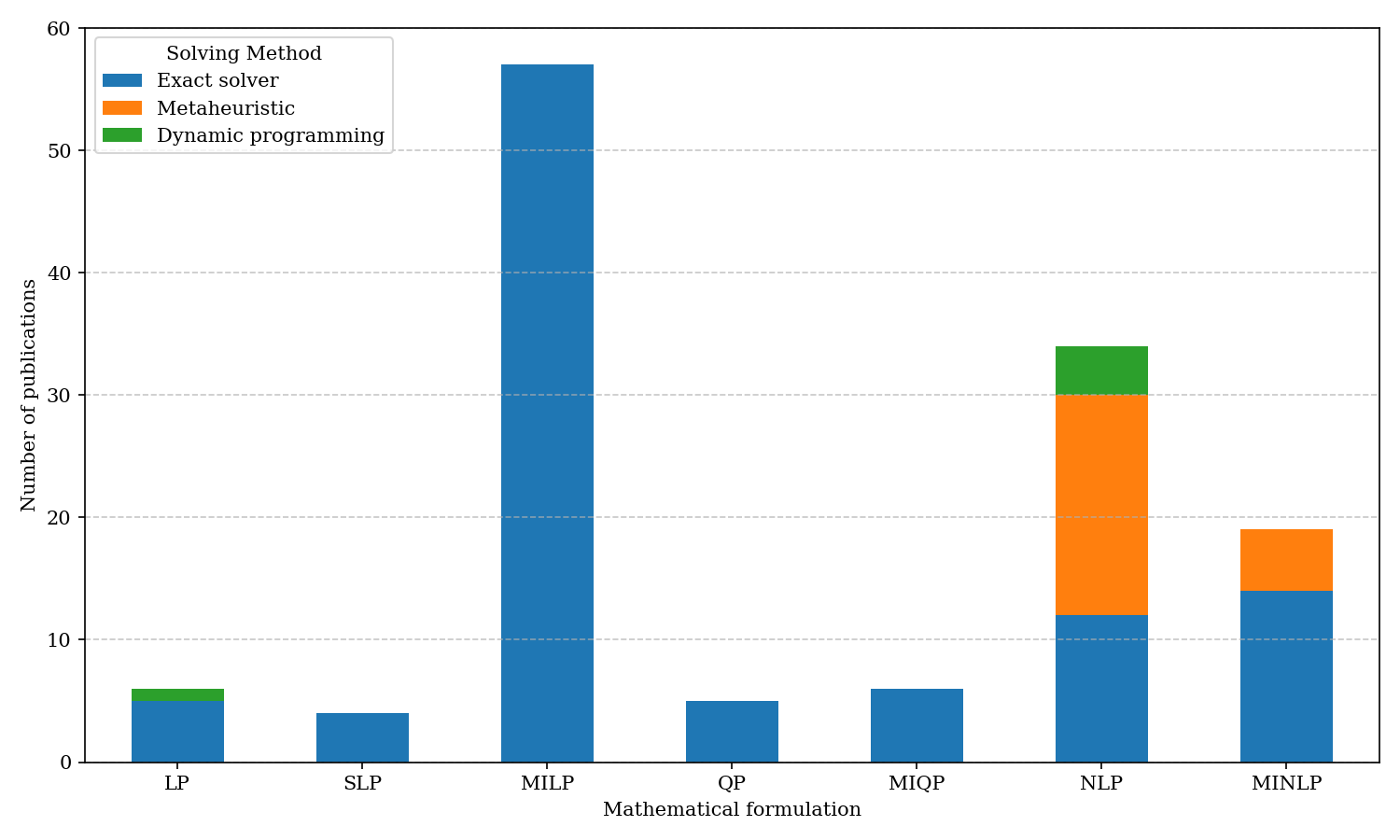} 
\caption{Distribution of mathematical modeling paradigms across the reviewed corpus.}
\label{fig:paradigm_distribution}
\end{figure}

Among the mathematical representations of hydropower cascading systems, \gls{milp} is the most prevalent \cite{8633891, 9211795, 7790824, 1137622, LIAO2024121341, FANG2024120932}, outpacing \gls{nlp} \cite{HERMIDA2018408, 4538514, Catalao2010Nonlinear, DoganMix, MPResOpt} and \gls{minlp} representations \cite{5565530, TeegavarapuSimonovic2000, HOSNAR2014194, SwarmMINLP, CATALAO2010904}, as shown in Fig.~\ref{fig:paradigm_distribution}. While \gls{minlp} and \gls{nlp} formulations allow for the most mathematically rigorous representation of nonlinear physical laws, they often rely on metaheuristic methods for resolution due to their extreme algebraic complexity \cite{SwarmMINLP, ZADEH20161393, NIU2021107315, RO-Xu, Wang2017OptOutputError}. This reliance frequently leads to locally optimal solutions rather than the provable global optima required for high-stakes interconnection studies.

Conversely, \gls{milp} has become the preferred bridge between accuracy and speed. Unlike standard \gls{lp} \cite{YOO2009182, MMLPNiu, MultiMohan, FENG2020119035, EFFLP, 8295134} or \gls{qp} \cite{7038354, 4682628, 7447809, 5953023, 6072098}, \gls{milp} utilizes binary variables to handle non-convex requirements such as \gls{huc} and turbine \gls{poz}. \gls{milp} formulations allow for the representation of \gls{pwl} functions, which serve as universal approximators of nonlinear physical relationships \cite{huang_relu_2020,vielma2010nonseparable, PLOUSSARD202450, ploussard2025tightening}. Furthermore, \gls{milp} formulations leverage mature solver algorithms that are generally faster and more reliable than nonlinear alternatives \cite{SU2020125556, ZHANG2019883,9887899,WANG2024130258}.

\subsection{The solver landscape: commercial vs. open-source}
\begin{figure}[ht!]
\centering
\includegraphics[width=1.0\textwidth]{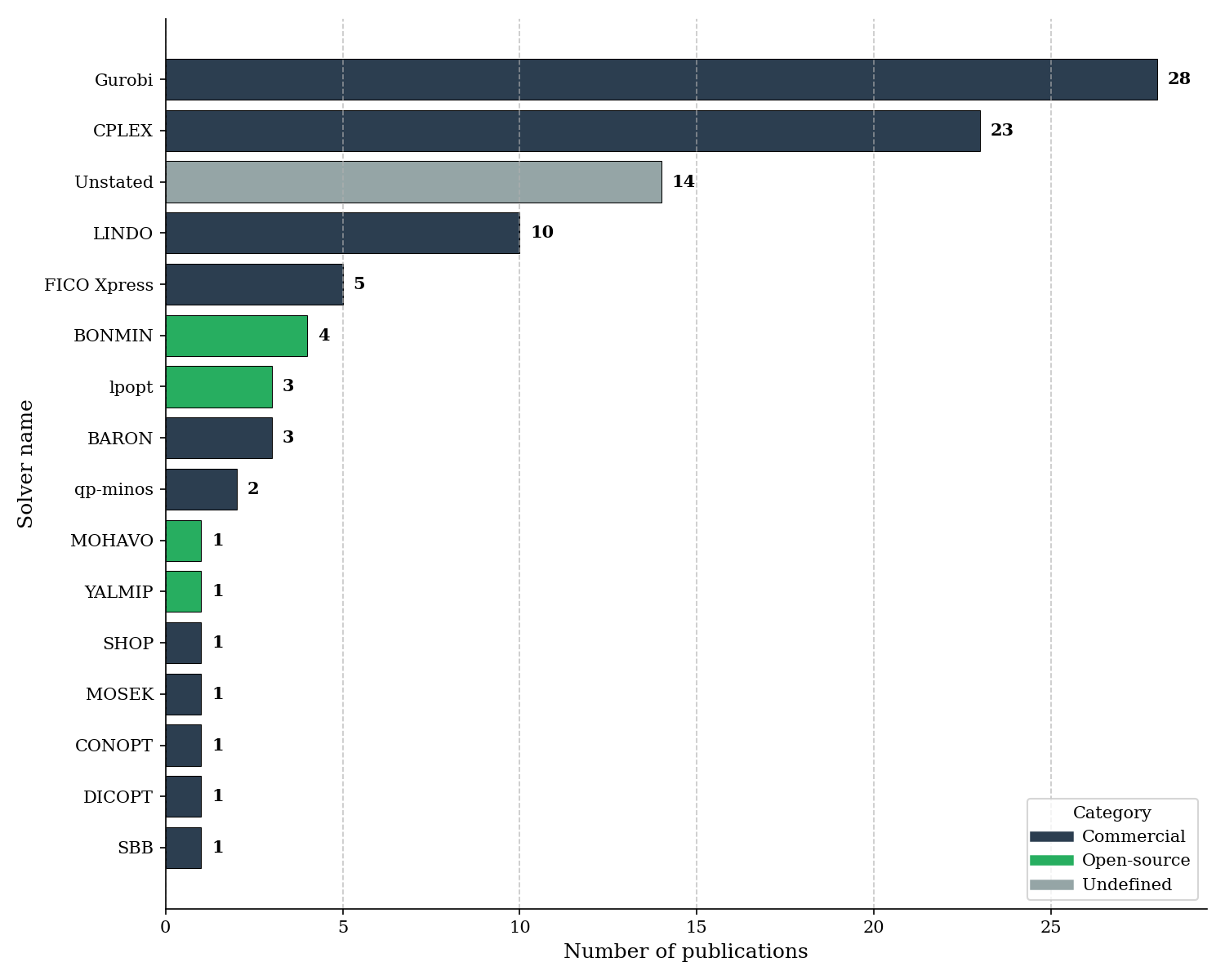} 
\caption{Distribution of numerical solvers categorized by commercial and open-source availability}
 \label{fig:solver_dist}
\end{figure}

The practical utility of these models is closely tied to the available software ecosystem. An analysis of the solver landscape in Fig.~\ref{fig:solver_dist} reveals that the vast majority of studies employ high-performance commercial solvers like CPLEX \cite{1137622, LU2024122085, 9314115, BELSNES2016167, 6575183} or Gurobi \cite{8633891, 9211795, 7790824, LIAO2024121341, FENG2022118620}. Conversely, the use of open-source solvers remains relatively limited in the literature \cite{Amina22, DoganMix, InflowInfluence, AMINABADI2024282, 7443262}, reflecting the comparatively lower algorithmic performance of many available alternatives, although recent advances have begun to narrow this gap. High reliance on proprietary tools, as opposed to open-source alternatives, may hinder broader adoption and slow the pace of model development and experimentation in this area due to licensing costs.

However, modern open-source solvers like SCIP \cite{Achterberg2009SCIP} or HiGHS \cite{Huangfu2018HiGHS}, which do not appear in the literature due to their recent maturity, have demonstrated increasingly competitive performance. These emerging tools provide a viable path for academic researchers to deploy sophisticated \gls{milp} and \gls{minlp} formulations while maintaining the provable optimality required by accurate hydropower cascading system models.

\section{Mathematical representation of hydropower cascading systems}

\label{sec:math}

This section establishes a standardized modeling framework, or ``Standard Model'', for cascading hydropower systems, defining a core set of ten fundamental physical and operational equations. It details the constituent variables and parameters while characterizing the mathematical nature of each constraint, specifically identifying linear, nonlinear, and non-convex components. The overview further evaluates the necessary adaptations to transform these physics-based rules into \gls{lp}, \gls{milp}, or \gls{minlp} formulations, consolidated in a taxonomic summary table. Finally, the section provides a systematic census of the 131 reviewed articles, analyzing how these specific equations are represented, simplified, or omitted within current literature.

\subsection{A standardized modeling framework}

To evaluate the technical evolution across the reviewed literature, it is necessary to establish a generalized mathematical framework that captures the physical and operational rules of a cascading system. This ``Standard Model'' serves as the baseline for our classification, allowing for a structured comparison of how different authors simplify, linearize, or expand upon these core constraints to address specific modeling gaps in power system studies.

The following formulation represents the integrated optimization problem for a cascading hydropower system:

\begin{align}
\max \quad & \mathcal{F} \label{eq:obj}\\[6pt]
\text{s.t.} \quad 
& V_{r,t} - V_{r,t-1} = (I_{r,t} - Q_{r,t} - L_{r,t}) \Delta t, 
&& \forall r \in \mathcal{R}, \, t \in \mathcal{T}, \label{eq:massbalance}\\[6pt]
& Q_{r,t} = \sum_{p \in \mathcal{P}(r)} Q^{{\text{P}}}_{p,t} + Q^{\text{NP}}_{r,t}, 
&& \forall r \in \mathcal{R}, \, t \in \mathcal{T}, \label{eq:totalrelease}\\[6pt]
& I_{r,t} = W_{r,t} + F_r\!\left(\{\, Q_{r',t'} : (r' \to r) \in \mathcal{A}, \, t' \leq t \,\}\right), 
&& \forall r \in \mathcal{R}, \, t \in \mathcal{T}, \label{eq:inflow}\\[6pt]
& E_{r,t} = f_r^{{\text{E}}}(V_{r,t-1},V_{r,t}), 
&& \forall r \in \mathcal{R}, \, t \in \mathcal{T}, \label{eq:elev-storage}\\[6pt]
& P_{p,t} = f_p^{{\text{P}}}(Q^{P}_{p,t}, h_{p,t}) 
&& \forall p \in \mathcal{P}, \, t \in \mathcal{T}, \label{eq:power}\\[6pt]
& h_{p,t} = E_{r(p),t} - f_r^{{\text{TR}}}(Q_{r(p),t}, E_{d(r(p)),t}) - f_p^{{\text{HL}}}(Q_{p,t}), 
&& \forall p \in \mathcal{P}, \, t \in \mathcal{T}, \label{eq:head}\\[6pt]
& E_{r,t} \in \left[\underline{E}_r, \overline{E}_r\right], 
&& \forall r \in \mathcal{R}, \, t \in \mathcal{T}, \label{eq:elevbounds}\\[6pt]
& (P_{p,t}, h_{p,t}) \in \mathcal{Z}_p, && \forall p \in \mathcal{P}, t \in \mathcal{T}, \label{eq:releaselimits}\\[6pt]
& Q^{\text{NP}}_{r,t} \geq 0, 
&& \forall r \in \mathcal{R}, \, t \in \mathcal{T}, \label{eq:nonpower}
\end{align}

The generalized objective function $\mathcal{F}$ in Equation~(\ref{eq:obj}) serves as the primary driver for the system's optimization. As analyzed in Section~\ref{sec:bibliometrics}, the literature has historically concentrated on the maximization of total energy production ($\mathcal{F} = \sum_{p \in \mathcal{P}, t \in \mathcal{T}} P_{p,t}$) or cumulative economic value ($\mathcal{F} = \sum_{p \in \mathcal{P}, t \in \mathcal{T}} \pi_{p,t} P_{p,t}$, where $\pi_{p,t}$ is the energy price). More recently, however, there has been a significant shift toward objectives for load peak shaving and grid stability enhancement. From a mathematical perspective, these formulations remain largely linear in nature, though they may adopt a quadratic structure when bilinear market or flow-head terms are introduced.

Equation~(\ref{eq:massbalance}) represents the reservoir mass-balance or water continuity constraint, which serves as the foundation for the spatiotemporal coupling within the cascade. This linear equation ensures that the variation in water storage $V_{r,t}$ over the interval $\Delta t$ is strictly governed by the conservation of water mass/volume: the volume increases with the total water inflow $I_{r,t}$ and decreases with the total reservoir release $Q_{r,t}$ and any parasitic losses $L_{r,t}$, such as evaporation or seepage. While the mass-balance itself is universally treated as an exact linear relationship, its operational accuracy is highly dependent on the fidelity of the inflow term $I_{r,t}$.

Equation~(\ref{eq:totalrelease}) constitutes a linear constraint that defines the total reservoir outflow $Q_{r,t}$ as the sum of power-generating discharges $Q_{p,t}^{{\text{P}}}$ and non-power releases $Q_{r,t}^{{\text{NP}}}$. This formulation accounts for the productive flow through all associated generating units $p \in \mathcal{P}(r)$ alongside non-productive flows, such as spillage (also called ``overflow'') or bypass through environmental outlets. To enhance computational tractability, modelers frequently aggregate individual turbine units into a single equivalent plant representation.

Equation~(\ref{eq:inflow}) characterizes the spatiotemporal coupling of the cascade by defining the total inflow $I_{r,t}$ as the sum of unregulated local runoff $W_{r,t}$ and the routing function $F_r$ of upstream releases. This function governs the transfer of water from an upstream facility $r'$ to the downstream reservoir $r$ across the set of hydraulic arcs $\mathcal{A}$, \textit{i.e.}, from all reservoirs $r'$ directly connected upstream of $r$. The function $F_r$ is generally modeled in three different ways: instantaneous transfer, fixed time delays, and distributed attenuation. In instantaneous transfer, water arrives within the same time step it was released ($Q_{r',t}$), whereas fixed delays shift the arrival by a specific duration $\tau$ ($Q_{r',t-\tau}$), and distributed delays spread a single pulse across multiple periods to model physical attenuation ($D_{(r' \rightarrow r),t}\star Q_{r',t}$, where $D_{(r' \rightarrow r),t}$ is the water travel time distribution from $r'$ to $r$). While Equation~(\ref{eq:inflow}) is generally linear, the formulation becomes nonlinear if the travel time distribution is modeled as a function of the discharge itself.

Equation~(\ref{eq:elev-storage}) models the reservoir storage-to-elevation relationship, a function determined strictly by the reservoir's specific geometric structure. Because $V_{r,t}$ denotes the storage volume at the end of period $t$, while $E_{r,t}$ represents the average forebay elevation during that interval, the elevation is most precisely modeled as a function of both current and preceding volumes. In practice, the literature predominantly defines $E_{r,t}$ as a function of the average volume: $(V_{r,t} + V_{r,t-1})/2$. Due to the shape of most reservoir basins, this function is generally concave nonlinear, though it can be linearized with high fidelity for small-scale volume fluctuations. For larger-scale fluctuations, it is typically either approximated by a high-order polynomial, preserving the model's nonlinear character, or by a \gls{pwl} function, allowing the use of \gls{milp} formulation.

Equation~(\ref{eq:power}) defines the power output $P_{p,t}$ of a unit as a function $f_p^{{\text{P}}}$ of its power-generating discharge $Q^{{\text{P}}}_{p,t}$ and its net hydraulic head $h_{p,t}$. This physical relationship, known as the fundamental hydropower equation, is commonly expressed as $P_{p,t} = \rho\, g\, \eta_{p,t} \, Q^{{\text{P}}}_{p,t} \,h_{p,t}$, where $\rho$ is the density of water, $g$ is the gravitational acceleration, and $\eta_{p,t}$ is the turbine efficiency. However, $\rho$ and $g$ are static parameters, and $\eta_{p,t}$ is a function of $Q^{{\text{P}}}_{p,t}$ and $h_{p,t}$, effectively making the power output a function of the turbine discharge and net hydraulic head. The function $f_p^{{\text{P}}}$ is inherently non-linear and non-convex. To manage this complexity, the literature frequently adopts a quadratic approximation by assuming a constant unit efficiency, or a strictly linear representation by further assuming a constant head in scenarios with minor reservoir elevation fluctuations. For high-fidelity applications where these simplifications are insufficient, researchers typically use \gls{pwl} approximations or \gls{slp} to navigate the non-convexity while maintaining computational tractability for large-scale cascades.

Equation~(\ref{eq:head}) defines the net hydraulic head $h_{p,t}$ by subtracting the tailwater level and the head loss from the reservoir elevation $E_{r(p),t}$. The tailwater level is theoretically a function $f_r^{{\text{TR}}}$ of both the downstream reservoir elevation $E_{d(r(p)),t}$ and the total reservoir discharge $Q_{r(p),t}$, although it is most commonly defined as a function of the total reservoir discharge only ($f_r^{{\text{TR}}}(Q_{r(p),t})$) for reservoirs that are far enough apart. The head loss, which represents the energy reduction of the water moving through the penstock, is a function of the unit release $Q_{p,t}$. From a mathematical perspective, Equation~(\ref{eq:head}) remains linear if the tailrace and head loss effects are omitted or represented by a first-order approximation; however, they introduce significant non-linearity when modeled as high-order functions. 

Equation~(\ref{eq:elevbounds}) establishes the feasible operational range of the reservoir by enforcing lower and upper elevation bounds, $\underline{E}_r$ and $\overline{E}_r$. These linear constraints represent critical physical and regulatory thresholds, such as flood control limits, minimum environmental flows, and structural safety margins.

Equation~(\ref{eq:releaselimits}) defines the feasible operating zone $\mathcal{Z}_p$ of the power output and net hydraulic head of generating unit $p$. Because of the discrete requirements of \gls{huc} and the presence of \gls{poz}, also known as ``rough zones'' where vibration and cavitation occur, this domain is inherently non-convex and composed of disjoint sets. Accurately capturing these disjoint sets typically necessitates the use of disjunctive constraints within an \gls{milp} formulation, utilizing binary variables to represent the "on/off" status and valid operating zones. This constraint can be linearized by relaxing the discrete operating requirements into a single, convex domain, which effectively assumes the unit is constantly online and ignores the physical turbulence of the rough zones. 

Equation~(\ref{eq:nonpower}) defines the non-power reservoir release $Q^{\text{NP}}_{r,t}$ as a non-negative variable representing all outflows that bypass the generating units, including spillway discharges and environmental bypass flows. Within the optimization framework, this variable often serves as a physical slack that ensures the feasibility of the mass-balance equation (Eq.~\ref{eq:massbalance}) during extreme hydrological events; because most objective functions prioritize productive generation, the solver naturally minimizes these non-productive releases unless the reservoir's physical or regulatory storage limits are reached.

It is important to emphasize that the ten equations detailed above characterize the fundamental physical and hydraulic behavior of the cascading system. They do not encompass the regulatory and environmental constraints that frequently govern real-world operational dispatch. 
Requirements such as discharge ramping limits (designed to protect downstream ecosystems and prevent riverbank erosion), maximum allowable daily fluctuations in reservoir levels, or mandated target release volumes for irrigation, navigation, and flood control are highly site-specific and non-universal. In practice, these regulatory and environmental rules are layered onto the physical backbone of the model as additional model constraints, tailored to the specific legal and environmental framework of the river basin under study.

\subsection{Mathematical model adaptations}

As established, the standard hydropower cascading model is a fundamentally non-convex, nonlinear mathematical program. The discrete requirements of \gls{huc} and the presence of \gls{poz} in Eq.~(\ref{eq:releaselimits}) introduce non-convex domains that necessitate the use of binary variables. Furthermore, the functions in Eq.~(\ref{eq:elev-storage}), Eq.~(\ref{eq:power}), and Eq.~(\ref{eq:head}) introduce significant nonlinearities. Consequently, a high-fidelity representation of a cascading system is inherently a \gls{minlp} model.

Avoiding integer variables requires ignoring disjunctive constraints and assuming the feasible operating domains of generating units to be convex (Eq.~(\ref{eq:releaselimits})). This also prevents the use of \gls{pwl} approximations to model nonlinear functions, although Eq.~(\ref{eq:power}) can still be represented as a trilinear or bilinear equation, and Eq.~(\ref{eq:elev-storage}) can be captured via high-fidelity polynomial approximations. Such a strategy effectively transforms the model into an \gls{nlp} or \gls{qp} problem, which avoids the combinatorial complexity of integer programming but may still struggle with non-convexity and local optima.

Alternatively, the model can be adapted into a high-fidelity \gls{milp} model. This approach preserves the discrete logic of \gls{huc} and \gls{poz} while utilizing arbitrarily accurate representations of nonlinear functions. Indeed, \gls{milp} formulations allow for the representation of any one-dimensional or multi-dimensional \gls{pwl} function \cite{vielma2010nonseparable}; these are recognized as universal approximators, sharing a mathematical foundation with modern \gls{relu}-based neural networks \cite{huang_relu_2020}. This paradigm is often preferred in the literature for its balance between physical fidelity and the availability of robust, deterministic solvers \cite{w11071392,taktak2017overview}.

Turning the model into a \gls{lp} model is significantly more restrictive, as it requires simplifications that substantially alter model fidelity. This includes not only ignoring \gls{huc} and \gls{poz}, but also assuming linear water delays (Eq.~(\ref{eq:inflow})), minor reservoir volume fluctuations (Eq.~(\ref{eq:elev-storage})), constant unit efficiency and reservoir hydraulic head (Eq.~(\ref{eq:power})), and negligible or linear tailrace and head loss effects (Eq.~(\ref{eq:head})). In some cases, \gls{lp} models can represent \gls{pwl} functions when they are concave and the objective promotes maximization of the function value \cite{7038354,7447809,FLETEN20082656,EFFLP,DELADURANTAYE2009499}. While computationally efficient, \gls{lp} adaptations are generally reserved for long-term planning or large-scale screening studies where unit-specific operational precision is less critical.

A summary of these model adaptations for each mathematical formulation is provided in Table~\ref{tab:eq_taxonomy}.

\begin{table}[ht!]
\centering
\caption{Adaptations of the hydropower cascading model equations for each mathematical formulation}
\label{tab:eq_taxonomy}
\scriptsize
\begin{tabularx}{\textwidth}{l p{1.5cm} p{3.2cm} p{3.2cm} p{3.2cm}}
\toprule
\textbf{Eq.} & \textbf{Nature} & \textbf{LP \newline Representation} & \textbf{MILP \newline Representation} & \textbf{NLP/MINLP \newline Representation} \\
\midrule

(\ref{eq:obj}) & Linear or nonlinear & 
Linear cost/benefit functions. & 
Includes binary penalties (e.g., start-up costs). & 
Nonlinear/stochastic objective functions. \\ \addlinespace

(\ref{eq:massbalance}) & Linear & 
Exact & 
Exact & 
Exact \\ \addlinespace

(\ref{eq:totalrelease}) & Linear & 
Exact & 
Exact & 
Exact \\ \addlinespace

(\ref{eq:inflow}) & Generally linear & 
Instantaneous transfer (no delay), fixed delay, or discrete convolution delay & 
Discrete selection between multiple possible time lags or binary-switched routing regimes & 
Continuous flow-dependent routing where the delay is a function of discharge \\ \addlinespace

(\ref{eq:elev-storage}) & Nonlinear & 
Constant (fixed head) or linear & 
\Gls{pwl} approximation & 
High-order polynomial \\ \addlinespace

(\ref{eq:power}) & Nonlinear & 
Linearized by fixing head and efficiency & 
\gls{1d} or \gls{2d} \gls{pwl} approximation & 
Bilinear product or function of discharge and net head \\ \addlinespace

(\ref{eq:head}) & Nonlinear & 
Ignored, constant, or linear tailrace and head loss & 
\Gls{pwl} approximation of tailrace-discharge relationship & 
Nonlinear tailrace-discharge function \\ \addlinespace

(\ref{eq:elevbounds}) & Linear & 
Exact linear bounds & 
Exact linear bounds & 
Exact linear bounds \\ \addlinespace

(\ref{eq:releaselimits}) & Non-convex & 
Power output relaxed to a single continuous interval (ignore rough zone) & 
Binary variables to enforce unit commitment or rough zone exclusion & 
Binary variables to enforce unit commitment or rough zone exclusion \\ \addlinespace

(\ref{eq:nonpower}) & Linear & 
Exact & 
Exact & 
Exact \\
\bottomrule
\end{tabularx}
\end{table}

\subsection{Comparative analysis of physical approximations}

Building upon the standardized framework and the mathematical paradigms defined previously, this subsection provides a systematic analysis of the modeling approximations employed across the 131 reviewed articles. Each fundamental equation is examined to identify prevalent trends in parameter and variable representation, and the mathematical representation of nonlinear functions. By contrasting the physical requirements of the cascading system with the practical constraint representation of existing models, this synthesis highlights the diverse strategies used in the literature to navigate the trade-off between model fidelity and computational tractability.

\subsubsection{Water balance equation and water loss representation (Eq.~ (\ref{eq:massbalance}))}

The water balance equation, often referred to as the mass balance or water continuity equation, serves as the fundamental equation for any reservoir model, providing a dynamic accounting of storage based on hydraulic flux. The primary drivers of this relationship are the cumulative inflows (comprising both natural side-inflows and upstream releases) and the total outflows (discharged through powerhouses, bypass, or spillways). While physical losses such as evaporation and seepage represent additional exit paths, they are frequently considered negligible in the literature due to their relatively minor impact on operational volumes compared to turbine discharge, especially for short time horizon (\textit{e.g.}, a few days) \cite{8633891, 5565530, 9211795, 7790824, 1137622, 7038354}.

In comparison, less than one-fifth of the reviewed articles explicitly account for such water losses. Among these works, the vast majority represent evaporation as a simplified, time-dependent parameter \cite{YOO2009182, MPResOpt, RO-Xu, YIN2022114582, Tayebiyan2016OptCleanEnergy}, and some authors extend this parametric representation to include seepage (leakage) \cite{Wang2017OptOutputError, Chen2023CascadeHydropower, 5287657}. From a physical standpoint, evaporative loss is more accurately modeled as a function of the water surface area, a dependency addressed in \cite{InflowInfluence} by modeling the surface area as a nonlinear (polynomial) function of the storage volume. Other researchers simplify this function by modeling evaporation directly as a function of the storage volume, using either nonlinear formulations \cite{SwarmMINLP} or linearized storage-to-evaporation approximations \cite{NLPOpt,MultiMohan}. 

\subsubsection{Representation of non-power releases (Eq.~ (\ref{eq:totalrelease}))}

Approximately one-quarter of the surveyed literature neglects non-power-generating discharges within their cascading models \cite{HOSNAR2014194, YUAN2022124025, 8494460, WANG202268, LIAO2021970, SU2020125556}. These studies typically prioritize the maximization of energy production or revenue, objectives for which spillage is viewed as a net loss or a non-contributing variable. In contrast, roughly three-quarters of the articles explicitly account for some form of non-generating release. The predominant form is spillage \cite{5953023,7853419,6019203,6072098,7483057}, which is also variously characterized in the literature as ``overflow'' \cite{MPResOpt,NLPOpt,Sharifi2021FDBMSA} or ``abandoned water'' \cite{LU2024122085,10853341}. Environmental bypass flows are also represented, though they are incorporated significantly less frequently within the reviewed corpus \cite{YOO2009182,OverflowMILP,ManikkuwahandiHornberger2021,SU2025124360}.

\subsubsection{Water routing equation and representation of water delay (Eq.~ (\ref{eq:inflow}))}

Approximately half of the reviewed articles assume no water delay between cascading reservoirs ($Q_{r',t}$), prioritizing computational simplicity \cite{5565530,HOSNAR2014194,LU2021114055,WU2023507,BELSNES2016167}. This assumption is generally justified when the model’s time step is significantly larger than the physical travel time of the water, rendering the delay negligible. Conversely, the remaining half of the literature utilizes a fixed time delay ($\tau$), represented as $Q_{r',t-\tau}$ \cite{LIAO2024121341,FANG2024120932,YUAN2022124025,8494460,WANG202268}. Beyond physical accuracy, the choice to omit or simplify delays is often driven by solver performance. From a computational perspective, the inclusion of time lags increases the bandwidth of the constraint matrix by spreading nonzero coefficients farther from the diagonal. This enlarged bandwidth typically leads to greater fill-in during sparse matrix factorization, increasing the density of the factorized matrices and, consequently, the computational cost and solution time \cite{davis2006direct}. Only two of the analyzed works incorporate the higher-fidelity discrete convolution approach ($D_{(r' \rightarrow r),t} \star Q_{r',t}$), which characterizes the spatiotemporal attenuation of the water pulse \cite{LIAO2021970,SLP01}. While this method is also linear and provides a more rigorous representation of hydraulic surges, it dramatically increases the number of non-zero entries in the model matrix, further compromising sparsity and significantly extending solution times.

\subsubsection{Storage-to-elevation equation (Eq.~(\ref{eq:elev-storage}))}

Strategies for modeling the reservoir storage-to-elevation relationship range from total omission to high-fidelity nonlinear representations. Several studies bypass the explicit volume-elevation link by using surrogate functions that represent either the hydraulic head or the total power output directly as a function of storage volume \cite{CATALAO2010904, 4682628, LIAO2021970, 9314115, 6039225}. To enhance computational tractability, some articles assume constant water levels for systems with low fluctuations \cite{7790824, 1137622, en12091604, AmaniHUC, MMLPNiu}, while others adopt a linear storage-elevation relationship as a compromise between accuracy and simplicity \cite{5565530, 7038354, FENG2022118620, YUAN2022124025, LU2021114055}. Many \gls{milp} formulations use \gls{pwl} representations enabled by \gls{sos2} functions or binary variables to capture the nonlinear nature of the relationship without sacrificing model linearity \cite{8633891, 9211795, LIAO2024121341, FANG2024120932, 8494460}. Finally, \gls{nlp} and \gls{minlp} models typically represent this relationship through high-order polynomials to maintain maximum physical fidelity \cite{TeegavarapuSimonovic2000, SwarmMINLP, Amina22, 6575183, ZADEH20161393}.

\subsubsection{Hydropower equation (Eq.~(\ref{eq:power}))}

The representation of the hydropower equation across the reviewed literature spans a wide spectrum of mathematical complexity, reflecting the inherent trade-off between physical fidelity and computational tractability. One simplification assumes both a constant hydraulic head and a constant unit efficiency, reducing the power output to a strictly linear function of the turbine discharge alone \cite{FENG2022118620, YUAN2022124025, LU2021114055, LU2024122085, BELSNES2016167}. While computationally attractive, this approximation is only physically justified when reservoir elevation fluctuations are negligible over the optimization horizon. 

A refinement of this approach retains the constant head and efficiency assumptions but represents the power-discharge relationship as a \gls{1d} \gls{pwl} function, capturing the nonlinear efficiency characteristics of the turbine across its operating range while preserving a linear or \gls{milp} structure \cite{8633891, 9211795, 7790824, 1137622, LIAO2024121341}. When the head is assumed to vary while the unit efficiency is still treated as fixed, the power output becomes a bilinear product of discharge and head, yielding a \gls{qp} formulation that introduces non-convexity and renders the problem harder to solve to global optimality \cite{HOSNAR2014194, SwarmMINLP, CATALAO2010904, 4682628, ZADEH20161393}. The most physically accurate representation models the power output as a fully nonlinear function of both discharge and net hydraulic head, accounting for the coupled efficiency surface of the turbine; however, this formulation requires \gls{nlp} or \gls{minlp} solvers that cannot guarantee global optimality for large-scale cascading instances \cite{5565530, TeegavarapuSimonovic2000, Amina22, 6575183, HERMIDA2018408}. 

To recover both physical fidelity and solver reliability, some articles employ \gls{2d} \gls{pwl} approximations of the power output surface as a function of both discharge and head, which can be embedded within a \gls{milp} formulation through the use of binary variables or \gls{sos2} \cite{7038354, FANG2024120932, WANG202268, LIAO2021970, SU2020125556}. Across the literature, the \gls{2d} \gls{pwl} approximation is consistently constructed using a rudimentary grid-based approach: each dimension is independently discretized into uniform intervals, forming a meshed grid that is subsequently triangulated
to define the affine pieces of the approximation. A less common and computationally lighter alternative is to represent the hydropower function by its convex hull, or McCormick envelope, expressed as a set of linear inequality constraints without any binary variables, yielding a pure \gls{lp} formulation \cite{10512930, su152416916, ZHANG2019883, SU2025124360, 8295134}. While more tractable, this relaxation introduces significant fidelity loss, as the convex hull simultaneously over- and underestimates the true feasible region of the flow-head-power relationship.

\subsubsection{Head loss and tailrace effects (Eq.~(\ref{eq:head}))}

The most prevalent approach in the surveyed literature is to omit both head loss and tailrace effects, simplifying the net hydraulic head to the difference between upstream and downstream reservoir elevations \cite{5565530, 1137622, HOSNAR2014194, CATALAO2010904, LU2021114055, 4682628}. Alternatively, some studies characterize head loss as a constant value \cite{7790824, 8494460, WANG202268, SU2020125556, WU2023507} or as a linear \cite{8633891, YUAN2022124025, 6575183, ZHAO2024129834, SU2025124360}, \gls{pwl} \cite{9211795, LIAO2024121341, FANG2024120932, LIAO2021970, SU2022395}, or nonlinear \cite{SwarmMINLP, Amina22, AMINABADI2024282, Liu2025DynamicPS, en18040964} function of the individual unit discharge. Similarly, tailrace effects are represented with varying degrees of complexity, ranging from constant parameters \cite{7790824, AmaniHUC, YIN2022114582, CHENG2022123908, FENG2020119035} to linear \cite{8633891, 7038354, FENG2022118620, YUAN2022124025, WANG202268}, \gls{pwl} \cite{9211795, LIAO2024121341, FANG2024120932, 8494460, LU2024122085}, or nonlinear \cite{TeegavarapuSimonovic2000, SwarmMINLP, Amina22, 6575183, ZADEH20161393} functions of the total reservoir discharge.

\subsubsection{HUC and POZ (Eq.~(\ref{eq:releaselimits}))}

The majority of the surveyed articles (approximately 60\%) omit both \gls{huc} and \gls{poz} constraints, generally favoring continuous formulations to ensure computational tractability \cite{LU2021126388,LIU2023127298,su151310002,LU202425,Xu2017multiobjective,KANG2026124031}. This simplification is mandatory for continuous model classes, such as \gls{lp}, \gls{qp}, and \gls{nlp}, which lack the binary variables required to represent these discrete operational thresholds. The remaining literature is evenly divided: roughly 20\% of the studies model \gls{huc} in isolation \cite{PlantsVsUnits,AMINABADI2024282,11234414,11009646,6009930,7804914}, while another 20\% incorporate both \gls{huc} and \gls{poz} to capture discrete operational limitations \cite{en14040887,en12091604,LUO2024110226,WU2024121502,10353562,su152416916}. For the vast majority of the articles modeling \gls{poz}, operating zones are represented by a set of disjoint feasible discharge or power output intervals. A select few authors \cite{9211795,SU2020125556,WU2023507} employ more advanced \gls{2d} operating zones, defining the feasible operating regions over both power output and net hydraulic head, capturing the interdependence between these two variables as a bounded area in the $(P_{p,t}, h_{p,t})$ plane. This approach captures the head-dependent nature of mechanical instabilities, such as cavitation and resonance, which are simplified in \gls{1d} representations. To integrate these often non-convex regions into a \gls{milp} framework, the \gls{2d} zones are partitioned into a set of convex domains (typically triangles). Each domain is then activated via binary variables, ensuring the unit operates strictly within safe hydraulic limits across the entire operational range of the reservoir.

A comprehensive summary of the modeling assumptions employed in each reviewed
article is provided in Table~\ref{tab:long_data}.

{\tiny
\begin{longtable}{p{0.5cm} p{1.6cm} p{0.5cm} p{0.5cm} p{1.5cm} p{0.3cm} p{0.4cm} p{0.6cm} p{0.7cm} p{1.1cm} p{1.0cm} p{0.5cm} p{0.5cm}}
\caption{Systematic summary of models} \label{tab:long_data} \\
\toprule
Article & Objective & Model & Type & Solver & Water Loss & Non-power & Water delay & Storage-to-elevation & Hydropower function & POZ & Head loss & Tail race \\
\midrule
\endfirsthead
\caption[]{Systematic summary of models} \\
\toprule
Article & Objective & Model & Type & Solver & Water loss & Non-power & Water delay & Storage-to-elevation & Hydropower function & POZ & Head loss & Tail race \\
\midrule
\endhead
\midrule
\multicolumn{13}{r}{Continued on next page} \\
\midrule
\endfoot
\bottomrule
\endlastfoot
\cite{8633891} & Economic & MILP & ES & Gurobi &  & \checkmark & Fixed & PWL & PWL & POZ & L & L \\
\cite{5565530} & Economic & MINLP & ES & SBB &  & \checkmark & No & L & NL &  &  &  \\
\cite{9211795} & Peak shaving & MILP & ES & Gurobi &  & \checkmark & No & PWL & PWL & 2D POZ & PWL & PWL \\
\cite{7790824} & Energy & MILP & ES & Gurobi &  & \checkmark & Fixed & C & PWL & POZ & C & C \\
\cite{1137622} & Economic & MILP & ES & CPLEX &  & \checkmark & Fixed & C & PWL & HUC &  &  \\
\cite{7038354} & Energy & QP & ES & qp-minos &  & \checkmark & No & L & 2D PWL* &  &  & L \\
\cite{TeegavarapuSimonovic2000} & Economic & MINLP & ES & DICOPT &  & \checkmark & Fixed & NL & NL &  &  & NL \\
\cite{HOSNAR2014194} & Economic & MINLP & ES & BARON &  &  & No &  & Q &  &  &  \\
\cite{LIAO2024121341} & Economic & MILP & ES & Gurobi &  & \checkmark & Fixed & PWL & PWL & POZ & PWL & PWL \\
\cite{SwarmMINLP} & Economic & MINLP & MH & - & \checkmark & \checkmark & No & NL & Q & HUC & NL & NL \\
\cite{FANG2024120932} & Peak shaving & MILP & ES & Unstated &  & \checkmark & Fixed & PWL & 2D PWL &  & PWL & PWL \\
\cite{FENG2022118620} & Economic & MILP & ES & Gurobi &  & \checkmark & No & L & L & POZ &  & L \\
\cite{YUAN2022124025} & Economic & MILP & ES & LINDO &  &  & Fixed & L & L & POZ & L & L \\
\cite{8494460} & Peak shaving & MILP & ES & LINDO &  &  & Fixed & PWL & PWL & POZ & C & PWL \\
\cite{CATALAO2010904} & Economic & MINLP & ES & FICO Xpress &  & \checkmark & Fixed &  & Q & POZ &  &  \\
\cite{LU2021114055} & Economic & MILP & ES & LINDO &  & \checkmark & No & L & L &  &  &  \\
\cite{4682628} & Economic & QP & ES & FICO Xpress &  & \checkmark & No &  & Q &  &  &  \\
\cite{LU2024122085} & Grid stability & MILP & ES & CPLEX &  & \checkmark & No & PWL & L &  &  & PWL \\
\cite{WANG202268} & Grid stability & MILP & ES & LINDO &  &  & Fixed & L & 2D PWL & HUC & C & L \\
\cite{LIAO2021970} & Grid stability & MIQP & ES & Gurobi &  &  & Conv &  & 2D PWL &  & PWL & PWL \\
\cite{SU2020125556} & Economic & MILP & ES & Gurobi &  &  & No & PWL & 2D PWL & 2D POZ & C & PWL \\
\cite{9314115} & Energy & MILP & ES & CPLEX &  & \checkmark & Fixed &  & 2D PWL &  &  &  \\
\cite{7447809} & Energy & QP & ES & qp-minos &  & \checkmark & Fixed & L & 2D PWL* &  &  & L \\
\cite{WU2023507} & Energy & MILP & ES & LINDO &  &  & No & PWL & PWL & 2D POZ & C & PWL \\
\cite{BELSNES2016167} & Economic & SLP & ES & CPLEX &  & \checkmark & No & L & L &  &  &  \\
\cite{Amina22} & Economic & MINLP & ES & BONMIN &  & \checkmark & No & NL & NL & HUC & NL & NL \\
\cite{10512930} & Economic & MILP & ES & Unstated &  &  & No & PWL & 2D CH &  & C & PWL \\
\cite{6575183} & Economic & MINLP & ES & CPLEX &  & \checkmark & Fixed & NL & NL &  & L & NL \\
\cite{HERMIDA2018408} & Economic & NLP & ES & CONOPT &  & \checkmark & Fixed & PWL* & NL &  &  &  \\
\cite{6039225} & Economic & MILP & ES & CPLEX &  & \checkmark & No &  & PWL &  &  &  \\
\cite{ZADEH20161393} & Economic & MINLP & MH & - &  & \checkmark & No & NL & Q &  & C & NL \\
\cite{SU2022395} & Peak shaving & MILP & ES & Gurobi &  &  & No & PWL & 2D PWL & POZ & PWL & PWL \\
\cite{en14040887} & Peak shaving & MILP & ES & CPLEX &  &  &  &  &  & POZ &  &  \\
\cite{962421} & Economic & MILP & ES & CPLEX &  & \checkmark & Fixed &  & PWL & HUC &  &  \\
\cite{4538514} & Economic & NLP & ES & FICO Xpress &  & \checkmark & No & L & Q &  &  &  \\
\cite{YOO2009182} & Energy & LP & ES & Unstated & \checkmark & \checkmark & No & L & L &  &  &  \\
\cite{en12091604} & Economic & MIQP & ES & CPLEX &  & \checkmark & Fixed & C & PWL & POZ &  &  \\
\cite{ZHAO2024129834} & Peak shaving & MILP & ES & Gurobi &  & \checkmark & Fixed & L & 2D PWL &  & L & L \\
\cite{SLP01} & Economic & SLP & ES & CPLEX &  & \checkmark & Conv & L & L &  &  & L \\
\cite{Catalao2010Nonlinear} & Economic & NLP & ES & FICO Xpress &  & \checkmark & No &  & Q &  &  &  \\
\cite{AmaniHUC} & Energy & MILP & ES & MOSEK &  &  & No & C & L & HUC & C & C \\
\cite{MMLPNiu} & Peak shaving & LP & ES & Unstated &  &  & No & C & L &  &  &  \\
\cite{FLETEN20082656} & Economic & MILP & ES & CPLEX &  & \checkmark & Fixed & C & PWL* & HUC &  & L \\
\cite{LUO2024110226} & Peak shaving & MILP & ES & CPLEX &  &  & No & PWL & 2D PWL & POZ & C & PWL \\
\cite{CHEN2023129185} & Energy & MILP & ES & Gurobi &  &  & No & C & C &  &  &  \\
\cite{DoganMix} & Economic & NLP & ES & lpopt &  & \checkmark & No & NL & Q &  &  &  \\
\cite{PlantsVsUnits} & Economic & MILP & ES & CPLEX &  &  &  &  & PWL & HUC & C &  \\
\cite{MPResOpt} & Economic & NLP & ES & Unstated & \checkmark & \checkmark & No &  & Q &  &  &  \\
\cite{ZhaoMILP} & Energy & MILP & ES & Gurobi &  & \checkmark & No &  & 2D PWL &  & PWL & PWL \\
\cite{OverflowMILP} & Economic & MILP & ES & SHOP &  & \checkmark & Fixed & L & PWL & HUC &  &  \\
\cite{NLPOpt} & Energy & NLP & ES & LINDO & \checkmark & \checkmark & No & NL & Q &  &  &  \\
\cite{AnderssonJohansson2025} & Economic & NLP & ES & Unstated &  & \checkmark & Fixed & NL & Q &  &  & L \\
\cite{MultiMohan} & Energy & LP & ES & Unstated & \checkmark & \checkmark & No &  & L &  &  &  \\
\cite{InflowInfluence} & Energy & NLP & ES & lpopt & \checkmark & \checkmark & No & NL & Q &  & C & NL \\
\cite{NIU2021107315} & Peak shaving & NLP & MH & - &  & \checkmark & No & NL & NL &  &  & NL \\
\cite{WU2024121502} & Peak shaving & MILP & ES & Gurobi &  &  & Fixed & PWL &  & POZ & C & PWL \\
\cite{RO-Xu} & Economic & MINLP & MH & - & \checkmark & \checkmark & No & NL & NL &  &  & NL \\
\cite{9887899} & Economic & MILP & ES & CPLEX &  & \checkmark & Fixed & PWL & PWL & HUC & PWL & PWL \\
\cite{AMINABADI2024282} & Economic & MINLP & ES & lpopt &  & \checkmark & No & NL & NL & HUC & NL & NL \\
\cite{10353562} & Peak shaving & MILP & ES & Unstated &  & \checkmark & Fixed & PWL & 2D PWL & POZ & PWL &  \\
\cite{DOGANIS20142051} & Economic & MIQP & ES & CPLEX &  &  & Fixed &  &  &  &  &  \\
\cite{Forknall2014} & Environmental & MILP & ES & LINDO &  & \checkmark & No &  &  &  &  &  \\
\cite{YIN2022114582} & Energy & NLP & ES & LINDO & \checkmark & \checkmark & No & NL & NL &  &  & C \\
\cite{Wang2017OptOutputError} & Energy & NLP & MH & - & \checkmark &  & No &  & NL &  &  &  \\
\cite{Liu2025DynamicPS} & Peak shaving & NLP & DP & - &  &  & Fixed & NL & NL &  & NL & NL \\
\cite{Zhu2025PFMODO} & Environmental & NLP & MH & - &  & \checkmark & No &  &  &  &  &  \\
\cite{NIU2018562} & Energy & NLP & MH & - &  & \checkmark & No & NL & NL &  &  & NL \\
\cite{LIAO2023127685} & Grid stability & MILP & ES & LINDO &  & \checkmark & Fixed & PWL & 2D PWL &  & PWL & PWL \\
\cite{Liu2016_JointOptimization} & Grid stability & MINLP & MH & - &  & \checkmark & Fixed &  & NL & HUC &  &  \\
\cite{xiao2023iboa} & Energy & NLP & MH & - &  & \checkmark &  & NL & NL &  &  & NL \\
\cite{Tayebiyan2016OptCleanEnergy} & Energy & NLP & MH & - & \checkmark & \checkmark & No & NL & NL &  &  &  \\
\cite{ManikkuwahandiHornberger2021} & Energy & NLP & MH & - & \checkmark & \checkmark & No & NL & NL &  &  &  \\
\cite{Chen2023CascadeHydropower} & Energy & NLP & MH & - & \checkmark &  & No &  &  &  &  &  \\
\cite{Zhang2019Cascade} & Energy & NLP & MH & - &  &  & No &  & NL &  &  &  \\
\cite{Niu2021CSA} & Energy & NLP & MH & - &  & \checkmark & No & NL & NL &  &  & NL \\
\cite{Sharifi2021FDBMSA} & Energy & NLP & MH & - & \checkmark & \checkmark & No & NL & NL &  &  & NL \\
\cite{Li2025CascadeHydroIntervals} & Energy & MILP & ES & Gurobi &  &  & No & C & L &  &  &  \\
\cite{ShangNested} & Energy & NLP & MH & - &  &  & No & NL & NL &  &  &  \\
\cite{en18040964} & Energy & NLP & MH & - &  & \checkmark & No & NL & NL &  & NL & NL \\
\cite{JointOpt} & Economic & NLP & MH & - &  &  & No &  &  &  &  &  \\
\cite{su152416916} & Peak shaving & MILP & ES & Gurobi &  & \checkmark & Fixed & PWL & 2D CH & POZ &  & PWL \\
\cite{627890} & Economic & MINLP & ES & Custom &  &  & Fixed &  & NL & POZ &  &  \\
\cite{QIU2021124239} & Energy & NLP & MH & - &  &  & No &  & NL &  &  &  \\
\cite{11234414} & Peak shaving & MILP & ES & CPLEX &  & \checkmark & Fixed & L & L & HUC &  &  \\
\cite{10672978} & Economic & NLP & DP & - &  & \checkmark & No &  & Q &  &  &  \\
\cite{7853589} & Economic & MINLP & ES & BARON &  & \checkmark & Fixed &  & NL & POZ &  &  \\
\cite{11009646} & Economic & MILP & ES & Custom &  & \checkmark & Fixed & L & L & HUC & C &  \\
\cite{4839958} & Energy & NLP & MH & - &  &  & Fixed &  & NL &  &  &  \\
\cite{5287657} & Energy & NLP & MH & - & \checkmark &  & No &  & NL &  &  &  \\
\cite{4601846} & Energy & NLP & ES & Custom &  & \checkmark & No & NL & NL &  & NL & NL \\
\cite{964906} & Peak shaving & MILP & ES & CPLEX &  & \checkmark & No & PWL & 2D PWL &  &  &  \\
\cite{6009930} & Economic & MILP & ES & CPLEX &  & \checkmark & Fixed &  & 2D PWL & HUC &  &  \\
\cite{7443262} & Economic & MINLP & ES & BONMIN &  & \checkmark & Fixed &  & NL &  &  &  \\
\cite{7804914} & Economic & MINLP & MH & - &  & \checkmark & No &  &  & HUC &  &  \\
\cite{5953023} & Economic & QP & ES & Unstated &  & \checkmark & No & L & Q &  &  &  \\
\cite{7853419} & Economic & MINLP & ES & BONMIN &  & \checkmark & No &  & NL & POZ &  &  \\
\cite{6019203} & Economic & MIQP & ES & FICO Xpress &  & \checkmark & No & L & Q & HUC &  &  \\
\cite{6072098} & Economic & QP & ES & Unstated &  & \checkmark & Fixed & L & Q &  &  &  \\
\cite{7483057} & Economic & MINLP & ES & BARON &  & \checkmark & Fixed &  & NL &  &  &  \\
\cite{7838997} & Economic & MINLP & ES & BONMIN &  & \checkmark & Fixed &  & NL & POZ &  &  \\
\cite{LIU202311} & Peak shaving & MILP & ES & Unstated &  & \checkmark & Fixed & PWL & 2D PWL &  & PWL & PWL \\
\cite{JIN2024132283} & Peak shaving & MILP & ES & Gurobi & \checkmark & \checkmark & Fixed & PWL & 2D PWL &  & PWL & PWL \\
\cite{ZHANG2019883} & Energy & MILP & ES & Gurobi &  & \checkmark & Fixed & PWL & 2D CH &  &  & PWL \\
\cite{MA2023138602} & Peak shaving & MILP & ES & Gurobi &  & \checkmark & No & PWL & 2D PWL &  & PWL & PWL \\
\cite{10853341} & Energy & NLP & DP & - &  & \checkmark & No &  & NL &  &  &  \\
\cite{LU2021126388} & Energy & NLP & MH & - &  & \checkmark & No & NL & NL &  &  &  \\
\cite{LIU2023127298} & Peak shaving & NLP & DP & - &  & \checkmark & No &  & NL &  &  &  \\
\cite{su151310002} & Economic & NLP & ES & Unstated &  & \checkmark & Fixed &  & PWL &  &  &  \\
\cite{LU202425} & Economic & MILP & ES & Gurobi &  &  & Fixed & L & PWL &  &  &  \\
\cite{LU2022435} & Economic & MIQP & ES & Gurobi &  &  & Fixed &  & L & HUC &  &  \\
\cite{ZHANG2023108833} & Economic & MILP & ES & Gurobi &  & \checkmark & No &  & L & HUC &  &  \\
\cite{Knezevic2016234245} & Economic & MILP & ES & Unstated &  & \checkmark & Fixed & C & PWL & HUC &  &  \\
\cite{Xu2017multiobjective} & Economic & NLP & ES & LINDO & \checkmark & \checkmark & Fixed & NL & NL &  &  &  \\
\cite{KANG2026124031} & Economic & MIQP & ES & Gurobi &  & \checkmark & Fixed & PWL & 2D PWL &  & PWL & PWL \\
\cite{LIU2025116775} & Peak shaving & MILP & ES & CPLEX &  & \checkmark & Fixed & PWL & 2D PWL & POZ & PWL & PWL \\
\cite{WANG2024130258} & Grid stability & MILP & ES & Gurobi &  &  & Fixed & PWL & 2D PWL & POZ &  & PWL \\
\cite{HUANG2025132756} & Grid stability & MILP & ES & YALMIP &  & \checkmark & No & PWL & PWL & POZ &  & PWL \\
\cite{en18143745} & Peak shaving & MILP & ES & Gurobi &  & \checkmark & No &  & 2D PWL &  & C &  \\
\cite{w17101441} & Peak shaving & SLP & ES & Gurobi &  &  & Fixed & PWL & 2D PWL &  & PWL & PWL \\
\cite{LI2023103534} & Grid stability & MINLP & ES & MOHAVO &  & \checkmark & No & NL & NL &  & NL & NL \\
\cite{en17112734} & Energy & MILP & ES & CPLEX &  & \checkmark & Fixed & PWL & 2D PWL & POZ & PWL & PWL \\
\cite{SU2025124360} & Economic & MILP & ES & Gurobi & \checkmark & \checkmark & Fixed & L & 2D CH & POZ & L & L \\
\cite{Shen24} & Peak shaving & MILP & ES & Gurobi &  & \checkmark & Fixed & PWL & 2D PWL & HUC & PWL & PWL \\
\cite{Metamodeling} & Economic & MILP & ES & CPLEX &  & \checkmark & Fixed & C & PWL &  &  &  \\
\cite{MaxConsume} & Energy & MILP & ES & CPLEX &  & \checkmark & Fixed & PWL & 2D PWL & POZ & PWL & PWL \\
\cite{CHENG2022123908} & Environmental & MILP & ES & Gurobi &  & \checkmark & Fixed & C & L & POZ &  & C \\
\cite{FENG2020119035} & Peak shaving & LP & DP & - &  & \checkmark & Fixed & C & L &  &  & C \\
\cite{EFFLP} & Energy & LP & ES & Gurobi &  & \checkmark & No &  & 2D PWL* &  &  & PWL \\
\cite{8295134} & Energy & LP & ES & Unstated &  & \checkmark & Fixed & PWL* & 2D CH &  &  &  \\
\cite{DELADURANTAYE2009499} & Economic & MILP & ES & CPLEX &  & \checkmark & Fixed & C & PWL* & HUC &  &  \\
\cite{LargeScaleOpt} & Energy & SLP & ES & MINOS &  & \checkmark & No & C & L &  &  & L \\ \midrule
\multicolumn{13}{p{\linewidth}}{\textit{ES: Exact Solver, MH: Metaheuristic, DP: Dynamic Programming \newline 
Conv: Convolutional \newline
C: Constant, L: Linear, PWL: piecewise linear, Q: Quadratic, NL: Nonlinear \newline
PWL*: concave piecewise linear function modeled using linear inequalities (no SOS2 or binary variables)}} \\
\bottomrule
\end{longtable}
}

\section{Discussion}

\label{sec:discussion}

\subsection{Balancing physical fidelity and computational complexity}

The selection of a modeling framework is governed by a fundamental trade-off between physical fidelity and computational tractability, a balance typically determined by the operational horizon of the study. High-fidelity components, such as tailrace effects (Eq.~\ref{eq:head}) to mitigate optimistic power capacity estimates and detailed routing (Eq.~\ref{eq:inflow}) to prevent inaccurate storage levels, are essential for short-term dispatch and reliability studies where precise hydraulic coupling is critical. Similarly, while the inclusion of discrete constraints like \gls{poz} (Eq.~\ref{eq:releaselimits}) introduces significant combinatorial overhead, it remains the only way to avoid release schedules that risk mechanical damage. While simplified linear models might still remain necessary for large systems and long-term planning where seasonal trends outweigh unit-specific nuances, the modern ``modeling gaps'' identified in this review highlight a growing necessity for high-fidelity representations to ensure that mathematical optimality aligns with physical reality in high-stakes environments.

\subsection{The dominance of MILP models and PWL representations}

A significant portion of the reviewed literature utilizes \gls{milp} to model hydropower cascading systems. This preference is primarily driven by the maturity and robustness of modern \gls{milp} solvers \cite{junger2010fifty}, such as Gurobi and CPLEX, which offer deterministic performance and guaranteed global optimality, features that are often elusive in general \gls{minlp} or non-convex \gls{nlp} solvers.

The flexibility of \gls{milp} stems from its ability to represent any \gls{pwl} function, either via \gls{sos2} for univariate cases \cite{williams_model_2013}, or through a combination of binary and continuous variables for multivariate functions \cite{vielma2010nonseparable}. Recognized as ``universal approximators'', \gls{pwl} functions can approximate nonlinear behavior at an arbitrary fidelity level \cite{huang_relu_2020}. This property bridges traditional optimization with modern machine learning architectures, specifically \gls{relu}-based neural networks. To ensure computational efficiency, ``tight'' \gls{pwl} formulations have been developed that significantly accelerate solver convergence \cite{vielma2010nonseparable}.

The fidelity and complexity of a \gls{pwl} representation are highly dependent on the approximation or regression method. Historically, \gls{pwl} regression relied on domain triangulation, \textit{i.e.} by partitioning the domain into uniformly distributed simplices \cite{7038354, FANG2024120932, WANG202268, LIAO2021970, SU2020125556}. More recently, however, two distinct paradigms have emerged in the literature that eliminate the need for predefined partitions, offering more adaptive and efficient approaches to \gls{pwl} approximation:
\begin{itemize}
\item The \gls{relu} neural network approach \cite{huang_relu_2020}: This method generates a \gls{pwl} approximation by training a \gls{relu} neural network on the target dataset. While computationally efficient and easily integrated into a \gls{milp} model \cite{GRIMSTAD2019106580}, it typically yields suboptimal approximations characterized by a relatively large number of affine pieces.
\item The \gls{milp}-based approach \cite{KAZDA2021107310,ploussard2025tightening}: This method identifies the optimal \gls{pwl} approximation of a dataset by solving a dedicated \gls{milp} problem. It generates the \gls{pwl} approximation with the minimum number of affine pieces required to satisfy a target error threshold. While mathematically rigorous, it can be computationally intensive, particularly for high-dimensional datasets \cite{PLOUSSARD202450, ploussard2025tightening}.
\end{itemize}

Both methods represent a significant methodological advancement over the
meshed grid representations consistently found in the cascading hydropower
literature. Whereas the grid-based approach constrains breakpoint placement
to a structured, axis-aligned partition of the discharge and head domains,
these two paradigms identify \gls{pwl} functions whose affine pieces are
defined over arbitrary domain partitions, enabling approximation effort to
be concentrated where the underlying nonlinear function exhibits the greatest
curvature. Furthermore, the \gls{milp}-based approach offers an additional
advantage of particular relevance to hydropower cascading models: by
minimizing the number of affine pieces required to meet a prescribed error
tolerance, it directly mitigates the curse of dimensionality associated with
the use of binary variables in \gls{milp} formulations, where each additional
affine piece introduces additional binary variables and tightening constraints
into the model.

\subsection{The rise of high-performance open-source solvers}

Historically, the literature on hydropower cascading optimization has been dominated by commercial solvers, most notably Gurobi and CPLEX. This trend is a direct consequence of the inherent complexity of cascading systems; as established in Section \ref{sec:math}, the interplay of nonlinear hydraulics and discrete operational constraints creates a high-dimensional, non-convex search space. Even when these physical relationships are simplified into \gls{lp} or \gls{milp} formulations via \gls{pwl} approximations, the resulting models often retain a significant computational overhead that, until recently, necessitated the specialized branching and cutting-edge heuristics found only in mature commercial solvers.

However, a notable gap exists in the current research landscape regarding the adoption of modern open-source alternatives. While tools like SCIP (for \gls{minlp} models) and HiGHS (for \gls{lp} and \gls{milp} models) have demonstrated competitive performance benchmarks that rival their commercial counterparts \cite{Achterberg2009SCIP,Huangfu2018HiGHS}, they remain largely absent from the 131 reviewed articles. This lack of representation suggests a significant opportunity for the field. A promising path for future research involves the integration of these powerful, freely available solvers into cascading models. By shifting away from a reliance on  proprietary licenses, the community can effectively lower the barrier to entry, fostering more reproducible research and allowing a broader range of institutions to tackle high-fidelity hydropower scheduling challenges.

\section{Conclusion}

\label{sec:conclusion}

This review has provided a comprehensive, equation-level analysis of 131 articles on rule-based mathematical programming for cascading hydropower systems, establishing a standardized 10-equation framework as a structured reference for classifying \gls{lp}, \gls{milp}, and \gls{minlp} formulations across the literature. By evaluating each article against every fundamental physical constraint (from water balance and spatiotemporal routing to head-dependent power generation and \gls{poz}) this census offers an empirical foundation for quantifying the current state of modeling fidelity in the field, a level of algebraic specificity absent from prior reviews.

A recurring theme across the reviewed corpus is the systematic omission or simplification of physical constraints that are critical for power grid reliability. Approximately 60\% of articles omit both \gls{poz} and \gls{huc}, while nominal hydraulic head assumptions and the neglect of tailrace and head loss effects remain pervasive. These simplifications, often justified by computational tractability, introduce a fidelity gap that this review frames explicitly through the lens of grid stability, a perspective aligned with the DOE HydroWIRES initiative and largely absent from prior works that focus primarily on water management performance.

The systematic census confirms that \gls{milp} is already the dominant paradigm in the literature, and this review reinforces its position as the optimal balance between physical accuracy and computational tractability. \gls{milp} formulations uniquely accommodate the non-convex discrete logic of \gls{huc} and \gls{poz} via binary variables, and the modeling of \gls{pwl} functions as universal approximators of the nonlinear hydraulic relationships. Critically, recent methodological advances have moved well beyond the rudimentary uniform grid-based triangulation that characterizes most \gls{pwl} representations in the cascading hydropower literature. Emerging approaches, including \gls{milp}-based optimal regression techniques, identify \gls{pwl} approximations with the minimum number of affine pieces needed to satisfy a prescribed error tolerance. By concentrating approximation effort where nonlinear curvature is greatest and directly minimizing the number of binary variables introduced, these methods directly mitigate the combinatorial overhead that has historically limited high-fidelity \gls{milp} models of cascading systems.

A bibliometric analysis of solver usage reveals that the literature is overwhelmingly dominated by commercial tools such as Gurobi and CPLEX, while high-performance open-source solvers (most notably SCIP and HiGHS) are virtually absent from the corpus despite their demonstrated competitive performance. This reliance on proprietary software creates a reproducibility barrier and limits access. The recent maturity of these open-source alternatives represents a critical and largely untapped opportunity: by enabling broader adoption of sophisticated \gls{milp} and \gls{minlp} formulations without licensing constraints, they can accelerate the development of standardized, high-fidelity cascading hydropower optimization frameworks that support both academic reproducibility and operational energy resilience.

\section*{Acknowledgements}

This work was authored for the Department of Energy (DOE) Office of Critical Minerals and Energy Innovation by Argonne National Laboratory, operated by UChicago Argonne, LLC, under contract number DE-AC02-06CH11357, and by Pacific Northwest National Laboratory, operated by Battelle Memorial Institute, under contract number DE-AC05-76RL01830. 
This study was supported by the HydroWIRES Initiative of DOE’s Hydropower and Hydrokinetics Office.

\section*{Declaration of generative AI and AI-assisted technologies in the manuscript preparation process}

During the preparation of this work, the authors used generative AI tools, including OpenAI's ChatGPT and Google's Gemini, to assist in refining the wording and improving the grammatical quality of the manuscript. After using these tools, the authors carefully reviewed and edited the content as needed and take full responsibility for the content of the published article.

\bibliographystyle{elsarticle-num}
\bibliography{cas-refs}

@inproceedings{11225846,
  author    = {Kincic, Slaven and Datta, Sohom and Nekkalapu, Sameer and Wang, Dewei and Mitra, Bhaskar and Vaiman, Marianna},
  booktitle = {2025 IEEE Power \& Energy Society General Meeting ({PESGM})}, 
  title     = {Potential Impact of Hydrological Modeling on Power System Reliability}, 
  year      = {2025},
  pages     = {1--5},
  doi       = {10.1109/PESGM52009.2025.11225846}
}

@inproceedings{11225344,
  author={Nekkalapu, Sameer and Datta, Sohom and Kincic, Slaven and Wang, Dewei and Mitra, Bhaskar and Mathew, Bibi},
  booktitle={2025 IEEE Power \& Energy Society General Meeting (PESGM)}, 
  title={Hydrological Modeling Impact on Western Interconnection Frequency Response Studies}, 
  year={2025},
  volume={},
  number={},
  pages={1-5},
  doi={10.1109/PESGM52009.2025.11225344}}

@misc{doe_reliability_topic,
  author       = {{U.S. Department of Energy}},
  title        = {Reliability},
  howpublished = "\url{https://www.energy.gov/topics/reliability}",
  year         = {2026}
}

@ARTICLE{8633891,
  author={Xia, Pei and Deng, Changhong and Chen, Yahong and Yao, Weiwei},
  journal={IEEE Access}, 
  title={{MILP} Based Robust Short-Term Scheduling for Wind–Thermal–Hydro Power System With Pumped Hydro Energy Storage}, 
  year={2019},
  volume={7},
  number={},
  pages={30261-30275},
  doi={10.1109/ACCESS.2019.2895090}}

@ARTICLE{5565530,
  author={Díaz, F. Javier and Contreras, Javier and Muñoz, José Ignacio and Pozo, David},
  journal={IEEE Transactions on Power Systems}, 
  title={Optimal Scheduling of a Price-Taker Cascaded Reservoir System in a Pool-Based Electricity Market}, 
  year={2011},
  volume={26},
  number={2},
  pages={604-615},
  doi={10.1109/TPWRS.2010.2063042}}

@ARTICLE{9211795,
  author={Zhao, Zhipeng and Cheng, Chuntian and Liao, Shengli and Li, Yapeng and Lü, Quan},
  journal={IEEE Transactions on Power Systems}, 
  title={A MILP Based Framework for the Hydro Unit Commitment Considering Irregular Forbidden Zone Related Constraints}, 
  year={2021},
  volume={36},
  number={3},
  pages={1819-1832},
  doi={10.1109/TPWRS.2020.3028480}}

@ARTICLE{7790824,
  author={Guedes, Lucas S. M. and de Mendonça Maia, Pedro and Lisboa, Adriano Chaves and Vieira, Douglas Alexandre Gomes and Saldanha, Rodney Rezende},
  journal={IEEE Transactions on Power Systems}, 
  title={A Unit Commitment Algorithm and a Compact MILP Model for Short-Term Hydro-Power Generation Scheduling}, 
  year={2017},
  volume={32},
  number={5},
  pages={3381-3390},
  doi={10.1109/TPWRS.2016.2641390}}

@ARTICLE{1137622,
  author={Conejo, A.J. and Arroyo, J.M. and Contreras, J. and Villamor, F.A.},
  journal={IEEE Transactions on Power Systems}, 
  title={Self-scheduling of a hydro producer in a pool-based electricity market}, 
  year={2002},
  volume={17},
  number={4},
  pages={1265-1272},
  doi={10.1109/TPWRS.2002.804951}}

@INPROCEEDINGS{7038354,
  author={Hamann, Andrew and Hug, Gabriela},
  booktitle={2014 Power Systems Computation Conference}, 
  title={Real-time optimization of a hydropower cascade using a linear modeling approach}, 
  year={2014},
  volume={},
  number={},
  pages={1-7},
  doi={10.1109/PSCC.2014.7038354}}

@article{TeegavarapuSimonovic2000,
  author    = {Teegavarapu, Ramesh S. V. and Simonovic, Slobodan P.},
  title     = {Short-term operation model for coupled hydropower reservoirs},
  journal   = {Journal of Water Resources Planning and Management},
  volume    = {126},
  number    = {2},
  pages     = {78--87},
  year      = {2000},
  publisher = {American Society of Civil Engineers (ASCE)},
  doi       = {10.1061/(ASCE)0733-9496(2000)126:2(78)},
}

@article{HOSNAR2014194,
title = {Mathematical modelling and MINLP programming of a hydro system for power generation},
journal = {Journal of Cleaner Production},
volume = {65},
pages = {194-201},
year = {2014},
issn = {0959-6526},
doi = {https://doi.org/10.1016/j.jclepro.2013.09.004},
url = {https://www.sciencedirect.com/science/article/pii/S0959652613005970},
author = {Jernej Hosnar and Anita Kovač-Kralj}
}

@article{LIAO2024121341,
title = {MILP model for short-term peak shaving of multi-grids using cascade hydropower considering discrete HVDC constraints},
journal = {Renewable Energy},
volume = {235},
pages = {121341},
year = {2024},
issn = {0960-1481},
doi = {https://doi.org/10.1016/j.renene.2024.121341},
url = {https://www.sciencedirect.com/science/article/pii/S0960148124014095},
author = {Shengli Liao and Jiang Xiong and Benxi Liu and Chuntian Cheng and Binbin Zhou and Yuqiang Wu}
}

@article{SwarmMINLP,
author = {Hatamkhani, Amir and Moridi, Ali and Yazdi, Jafar},
year = {2019},
month = {12},
pages = {},
title = {A Simulation – Optimization Models for Multi-Reservoir Hydropower Systems Design at Watershed Scale},
volume = {149},
journal = {Renewable Energy},
doi = {10.1016/j.renene.2019.12.055}
}

@article{FANG2024120932,
title = {An MILP model based on a processing strategy of complex multisource constraints for the short-term peak shaving operation of large-scale cascaded hydropower plants},
journal = {Renewable Energy},
volume = {231},
pages = {120932},
year = {2024},
issn = {0960-1481},
doi = {https://doi.org/10.1016/j.renene.2024.120932},
url = {https://www.sciencedirect.com/science/article/pii/S0960148124010000},
author = {Zhou Fang and Shengli Liao and Hongye Zhao and Chuntian Cheng and Benxi Liu and Huan Wang and Shushan Li}
}

@article{FENG2022118620,
title = {Weekly hydropower scheduling of cascaded reservoirs with hourly power and capacity balances},
journal = {Applied Energy},
volume = {311},
pages = {118620},
year = {2022},
issn = {0306-2619},
doi = {https://doi.org/10.1016/j.apenergy.2022.118620},
url = {https://www.sciencedirect.com/science/article/pii/S0306261922000940},
author = {Suzhen Feng and Hao Zheng and Yifan Qiao and Zetai Yang and Jinwen Wang and Shuangquan Liu}
}

@article{YUAN2022124025,
title = {Optimal scheduling of cascade hydropower plants in a portfolio electricity market considering the dynamic water delay},
journal = {Energy},
volume = {252},
pages = {124025},
year = {2022},
issn = {0360-5442},
doi = {https://doi.org/10.1016/j.energy.2022.124025},
url = {https://www.sciencedirect.com/science/article/pii/S0360544222009288},
author = {Wenlin Yuan and Shijie Zhang and Chengguo Su and Yang Wu and Denghua Yan and Zening Wu}
}

@INPROCEEDINGS{8494460,
  author={Su, Chengguo and Cheng, Chuntian and Wang, Peilin},
  booktitle={2018 IEEE International Conference on Environment and Electrical Engineering and 2018 IEEE Industrial and Commercial Power Systems Europe (EEEIC / I\&CPS Europe)}, 
  title={An MILP Model for Short-Term Peak Shaving Operation of Cascaded Hydropower Plants Considering Unit Commitment}, 
  year={2018},
  volume={},
  number={},
  pages={1-5},
  doi={10.1109/EEEIC.2018.8494460}}

@article{CATALAO2010904,
title = {Scheduling of head-dependent cascaded reservoirs considering discharge ramping constraints and start/stop of units},
journal = {International Journal of Electrical Power \& Energy Systems},
volume = {32},
number = {8},
pages = {904-910},
year = {2010},
issn = {0142-0615},
doi = {https://doi.org/10.1016/j.ijepes.2010.01.022},
url = {https://www.sciencedirect.com/science/article/pii/S0142061510000311},
author = {J.P.S. Catalão and H.M.I. Pousinho and V.M.F. Mendes}
}

@article{LU2021114055,
title = {Optimization model for the short-term joint operation of a grid-connected wind-photovoltaic-hydro hybrid energy system with cascade hydropower plants},
journal = {Energy Conversion and Management},
volume = {236},
pages = {114055},
year = {2021},
issn = {0196-8904},
doi = {https://doi.org/10.1016/j.enconman.2021.114055},
url = {https://www.sciencedirect.com/science/article/pii/S0196890421002314},
author = {Lu Lu and Wenlin Yuan and Chengguo Su and Peilin Wang and Chuntian Cheng and Denghua Yan and Zening Wu}
}

@ARTICLE{4682628,
  author={Catalao, J. P. S. and Mariano, S. J. P. S. and Mendes, V. M. F. and Ferreira, L. A. F. M.},
  journal={IEEE Transactions on Power Systems}, 
  title={Scheduling of Head-Sensitive Cascaded Hydro Systems: A Nonlinear Approach}, 
  year={2009},
  volume={24},
  number={1},
  pages={337-346},
  doi={10.1109/TPWRS.2008.2005708}}

@article{LU2024122085,
title = {Medium- and long-term interval optimal scheduling of cascade hydropower-photovoltaic complementary systems considering multiple uncertainties},
journal = {Applied Energy},
volume = {353},
pages = {122085},
year = {2024},
issn = {0306-2619},
doi = {https://doi.org/10.1016/j.apenergy.2023.122085},
url = {https://www.sciencedirect.com/science/article/pii/S0306261923014496},
author = {Na Lu and Guangyan Wang and Chengguo Su and Zaimin Ren and Xiaoyue Peng and Quan Sui}
}

@article{WANG202268,
title = {Short-term optimal scheduling of cascade hydropower plants shaving peak load for multiple power grids},
journal = {Renewable Energy},
volume = {184},
pages = {68-79},
year = {2022},
issn = {0960-1481},
doi = {https://doi.org/10.1016/j.renene.2021.10.079},
url = {https://www.sciencedirect.com/science/article/pii/S096014812101538X},
author = {Peilin Wang and Wenlin Yuan and Chengguo Su and Yang Wu and Lu Lu and Denghua Yan and Zening Wu}
}

@article{LIAO2021970,
title = {Daily peak shaving operation of cascade hydropower stations with sensitive hydraulic connections considering water delay time},
journal = {Renewable Energy},
volume = {169},
pages = {970-981},
year = {2021},
issn = {0960-1481},
doi = {https://doi.org/10.1016/j.renene.2021.01.072},
url = {https://www.sciencedirect.com/science/article/pii/S0960148121000793},
author = {Shengli Liao and Zhanwei Liu and Benxi Liu and Chuntian Cheng and Xinyu Wu and Zhipeng Zhao}
}

@article{SU2020125556,
title = {Short-term generation scheduling of cascade hydropower plants with strong hydraulic coupling and head-dependent prohibited operating zones},
journal = {Journal of Hydrology},
volume = {591},
pages = {125556},
year = {2020},
issn = {0022-1694},
doi = {https://doi.org/10.1016/j.jhydrol.2020.125556},
url = {https://www.sciencedirect.com/science/article/pii/S0022169420310167},
author = {Chengguo Su and Wenlin Yuan and Chuntian Cheng and Peilin Wang and Lifei Sun and Taiheng Zhang}
}

@ARTICLE{9314115,
  author={Xie, Jun and Zheng, Yimin and Pan, Xueping and Zheng, Yuan and Zhang, Liqin and Zhan, Yongsheng},
  journal={IEEE Access}, 
  title={A Short-Term Optimal Scheduling Model for Wind-Solar-Hydro Hybrid Generation System With Cascade Hydropower Considering Regulation Reserve and Spinning Reserve Requirements}, 
  year={2021},
  volume={9},
  number={},
  pages={10765-10777},
  doi={10.1109/ACCESS.2021.3049280}}

@ARTICLE{7447809,
  author={Hamann, Andrew and Hug, Gabriela and Rosinski, Stan},
  journal={IEEE Transactions on Power Systems}, 
  title={Real-Time Optimization of the Mid-Columbia Hydropower System}, 
  year={2017},
  volume={32},
  number={1},
  pages={157-165},
  doi={10.1109/TPWRS.2016.2550490}}

@article{WU2023507,
title = {A mixed-integer linear programming model for hydro unit commitment considering operation constraint priorities},
journal = {Renewable Energy},
volume = {204},
pages = {507-520},
year = {2023},
issn = {0960-1481},
doi = {https://doi.org/10.1016/j.renene.2023.01.006},
url = {https://www.sciencedirect.com/science/article/pii/S096014812300006X},
author = {Xinyu Wu and Yiyang Wu and Xilong Cheng and Chuntian Cheng and Zehong Li and Yongqi Wu}
}

@article{BELSNES2016167,
title = {Applying successive linear programming for stochastic short-term hydropower optimization},
journal = {Electric Power Systems Research},
volume = {130},
pages = {167-180},
year = {2016},
issn = {0378-7796},
doi = {https://doi.org/10.1016/j.epsr.2015.08.020},
url = {https://www.sciencedirect.com/science/article/pii/S0378779615002576},
author = {M.M. Belsnes and O. Wolfgang and T. Follestad and E.K. Aasgård}
}

@article{Amina22,
author = {Aminabadi, M. and Séguin, Sara and Fofana, Issouf and Fleten, S.-E and Aasgård, Ellen},
year = {2023},
month = {09},
pages = {1-25},
title = {Short-term hydropower optimization in the day-ahead market using a nonlinear stochastic programming model},
journal = {Energy Systems},
doi = {10.1007/s12667-023-00618-8}
}

@INPROCEEDINGS{10512930,
  author={Ge, Xiaohui and Zhang, Xuesong and Fang, Ni and Jiao, Ye and Zhao, Bo and Li, Yang and Wu, Feng},
  booktitle={2023 IEEE 7th Conference on Energy Internet and Energy System Integration (EI2)}, 
  title={Short-Term Scheduling of Cascaded Hydropower Stations with Pumped Storage Considering Strong Hydraulic Coupling}, 
  year={2023},
  volume={},
  number={},
  pages={3311-3316},
  doi={10.1109/EI259745.2023.10512930}}

@ARTICLE{6575183,
  author={Lima, Ricardo M. and Marcovecchio, Marian G. and Novais, Augusto Queiroz and Grossmann, Ignacio E.},
  journal={IEEE Transactions on Power Systems}, 
  title={On the Computational Studies of Deterministic Global Optimization of Head Dependent Short-Term Hydro Scheduling}, 
  year={2013},
  volume={28},
  number={4},
  pages={4336-4347},
  doi={10.1109/TPWRS.2013.2274559}}

@article{HERMIDA2018408,
title = {On the hydropower short-term scheduling of large basins, considering nonlinear programming, stochastic inflows and heavy ecological restrictions},
journal = {International Journal of Electrical Power \& Energy Systems},
volume = {97},
pages = {408-417},
year = {2018},
issn = {0142-0615},
doi = {https://doi.org/10.1016/j.ijepes.2017.10.033},
url = {https://www.sciencedirect.com/science/article/pii/S0142061517310323},
author = {Gloria Hermida and Edgardo D. Castronuovo}
}

@INPROCEEDINGS{6039225,
  author={Tong, Bo and Guan, Xiaohong and Zhai, Qiaozhu and Gao, Feng},
  booktitle={2011 IEEE Power and Energy Society General Meeting}, 
  title={Long-term scheduling of cascaded hydro energy system with distributed water usage allocation constraints}, 
  year={2011},
  volume={},
  number={},
  pages={1-7},
  doi={10.1109/PES.2011.6039225}}

@article{ZADEH20161393,
title = {Optimal Design and Operation of Hydraulically Coupled Hydropower Reservoirs System},
journal = {Procedia Engineering},
volume = {154},
pages = {1393-1400},
year = {2016},
note = {12th International Conference on Hydroinformatics (HIC 2016) - Smart Water for the Future},
issn = {1877-7058},
doi = {https://doi.org/10.1016/j.proeng.2016.07.509},
url = {https://www.sciencedirect.com/science/article/pii/S1877705816318987},
author = {N. Afsharian Zadeh and S.J. Mousavi and E. Jahani and J.H. Kim}
}

@article{SU2022395,
title = {Short-term optimal scheduling of cascade hydropower plants with reverse-regulating effects},
journal = {Renewable Energy},
volume = {199},
pages = {395-406},
year = {2022},
issn = {0960-1481},
doi = {https://doi.org/10.1016/j.renene.2022.08.159},
url = {https://www.sciencedirect.com/science/article/pii/S0960148122013362},
author = {Chengguo Su and Peilin Wang and Wenlin Yuan and Yang Wu and Feng Jiang and Zening Wu and Denghua Yan}
}

@Article{en14040887,
AUTHOR = {Cheng, Xianliang and Feng, Suzhen and Huang, Yanxuan and Wang, Jinwen},
TITLE = {A New Peak-Shaving Model Based on Mixed Integer Linear Programming with Variable Peak-Shaving Order},
JOURNAL = {Energies},
VOLUME = {14},
YEAR = {2021},
NUMBER = {4},
ARTICLE-NUMBER = {887},
URL = {https://www.mdpi.com/1996-1073/14/4/887},
ISSN = {1996-1073},
DOI = {10.3390/en14040887}
}

@ARTICLE{962421,
  author={Chang, G.W. and Aganagic, M. and Waight, J.G. and Medina, J. and Burton, T. and Reeves, S. and Christoforidis, M.},
  journal={IEEE Transactions on Power Systems}, 
  title={Experiences with mixed integer linear programming based approaches on short-term hydro scheduling}, 
  year={2001},
  volume={16},
  number={4},
  pages={743-749},
  doi={10.1109/59.962421}}

@INPROCEEDINGS{4538514,
  author={Mariano, S. J. P. S. and Catalao, J. P. S. and Mendes, V. M. F. and Ferreira, L. A. F. M.},
  booktitle={2007 IEEE Lausanne Power Tech}, 
  title={Profit-Based Short-Term Hydro Scheduling considering Head-Dependent Power Generation}, 
  year={2007},
  volume={},
  number={},
  pages={1362-1367},
  doi={10.1109/PCT.2007.4538514}}

@article{YOO2009182,
title = {Maximization of hydropower generation through the application of a linear programming model},
journal = {Journal of Hydrology},
volume = {376},
number = {1},
pages = {182-187},
year = {2009},
issn = {0022-1694},
doi = {https://doi.org/10.1016/j.jhydrol.2009.07.026},
url = {https://www.sciencedirect.com/science/article/pii/S0022169409004193},
author = {Ju-Hwan Yoo}
}

@Article{en12091604,
AUTHOR = {Fekete, Krešimir and Nikolovski, Srete and Klaić, Zvonimir and Androjić, Ana},
TITLE = {Optimal Re-Dispatching of Cascaded Hydropower Plants Using Quadratic Programming and Chance-Constrained Programming},
JOURNAL = {Energies},
VOLUME = {12},
YEAR = {2019},
NUMBER = {9},
ARTICLE-NUMBER = {1604},
URL = {https://www.mdpi.com/1996-1073/12/9/1604},
ISSN = {1996-1073},
DOI = {10.3390/en12091604}
}

@article{ZHAO2024129834,
title = {Short-term peak-shaving operation of “N-reservoirs and multicascade” large-scale hydropower systems based on a decomposition-iteration strategy},
journal = {Energy},
volume = {288},
pages = {129834},
year = {2024},
issn = {0360-5442},
doi = {https://doi.org/10.1016/j.energy.2023.129834},
url = {https://www.sciencedirect.com/science/article/pii/S0360544223032280},
author = {Hongye Zhao and Shengli Liao and Zhou Fang and Benxi Liu and Xiangyu Ma and Jia Lu}
}

@article{SLP01,
author = {Gauvin, Charles and Delage, Erick and Gendreau, Michel},
year = {2018},
month = {01},
pages = {},
title = {A successive linear programming algorithm with non-linear time series for the reservoir management problem},
volume = {15},
journal = {Computational Management Science},
doi = {10.1007/s10287-017-0295-4}
}

@article{Catalao2010Nonlinear, 
author = {Catalao, J. P. S. and Mariano, S. J. P. S. and Mendes, V. M. F.  and Ferreira, L. A. F. M. }, 
title = {Nonlinear optimization method for short-term hydro scheduling considering head-dependency}, 
journal = {European Transactions on Electrical Power}, 
year = {2010}, 
volume = {20}, 
pages = {172–183},
doi = {10.1002/etep.301}, 
note = {Published online 13 November 2008}, 
url = {https://onlinelibrary.wiley.com/doi/10.1002/etep.301} }

@article{AmaniHUC,
author = {Amani, Alireza and Alizadeh, Hosein},
year = {2021},
month = {04},
pages = {1-19},
title = {Solving Hydropower Unit Commitment Problem Using a Novel Sequential Mixed Integer Linear Programming Approach},
volume = {35},
journal = {Water Resources Management},
doi = {10.1007/s11269-021-02806-6}
}

@article{MMLPNiu,
author = {Wen-jing Niu  and Zhong-kai Feng  and Chun-tian Cheng },
title = {Min-Max Linear Programming Model for Multireservoir System Operation with Power Deficit Aspect},
journal = {Journal of Water Resources Planning and Management},
volume = {144},
number = {10},
pages = {06018006},
year = {2018},
doi = {10.1061/(ASCE)WR.1943-5452.0000977},

URL = {https://ascelibrary.org/doi/abs/10.1061/(ASCE)WR.1943-5452.0000977},
eprint = {https://ascelibrary.org/doi/pdf/10.1061/(ASCE)WR.1943-5452.0000977}
}

@article{FLETEN20082656,
title = {Short-term hydropower production planning by stochastic programming},
journal = {Computers \& Operations Research},
volume = {35},
number = {8},
pages = {2656-2671},
year = {2008},
note = {Queues in Practice},
issn = {0305-0548},
doi = {https://doi.org/10.1016/j.cor.2006.12.022},
url = {https://www.sciencedirect.com/science/article/pii/S0305054806003224},
author = {Stein-Erik Fleten and Trine Krogh Kristoffersen}
}

@article{LUO2024110226,
title = {Short-term peak shaving model of cascade hybrid pumped storage hydropower station retrofitted from conventional hydropower},
journal = {International Journal of Electrical Power \& Energy Systems},
volume = {162},
pages = {110226},
year = {2024},
issn = {0142-0615},
doi = {https://doi.org/10.1016/j.ijepes.2024.110226},
url = {https://www.sciencedirect.com/science/article/pii/S0142061524004472},
author = {Bin Luo and Xinyu Liu and Yongcan Chen and Can Zhou and Xin Long}
}

@article{CHEN2023129185,
title = {A stochastic linear programming model for maximizing generation and firm output at a reliability in long-term hydropower reservoir operation},
journal = {Journal of Hydrology},
volume = {618},
pages = {129185},
year = {2023},
issn = {0022-1694},
doi = {https://doi.org/10.1016/j.jhydrol.2023.129185},
url = {https://www.sciencedirect.com/science/article/pii/S0022169423001270},
author = {Cheng Chen and Suzhen Feng and Shuangquan Liu and Hao Zheng and Hong Zhang and Jinwen Wang}
}

@article{DoganMix,
author = {Mustafa S. Dogan  and Jay R. Lund  and Josue Medellin-Azuara },
title = {Hybrid Linear and Nonlinear Programming Model for Hydropower Reservoir Optimization},
journal = {Journal of Water Resources Planning and Management},
volume = {147},
number = {3},
pages = {06021001},
year = {2021},
doi = {10.1061/(ASCE)WR.1943-5452.0001353},

URL = {https://ascelibrary.org/doi/abs/10.1061/(ASCE)WR.1943-5452.0001353},
eprint = {https://ascelibrary.org/doi/pdf/10.1061/(ASCE)WR.1943-5452.0001353}
}

@article{PlantsVsUnits,
author = {Pérez-Díaz, Juan and González-Martínez, Jesús},
year = {2023},
month = {12},
pages = {1-21},
title = {Comparison between plant-based and unit-based production functions for the day-ahead energy and reserve scheduling of a hydropower plant with common waterways},
journal = {Energy Systems},
doi = {10.1007/s12667-023-00643-7}
}

@article{MPResOpt,
author = {Nguyen, Dung},
year = {2023},
month = {03},
pages = {5444},
title = {Operating Multi-Purpose Reservoirs in the Red River Basin: Hydropower Benefit Optimization in Conditions Ensuring Enough Water for Downstream Irrigation},
volume = {15},
journal = {Sustainability},
doi = {10.3390/su15065444}
}

@article{ZhaoMILP,
author = {Zhao, Zhipeng and Cheng, ChunTian and Yan, Lingzhi},
title = {An efficient and accurate Mixed-integer linear programming model for long-term operations of large-scale hydropower systems},
journal = {IET Renewable Power Generation},
volume = {15},
number = {6},
pages = {1178-1190},
doi = {https://doi.org/10.1049/rpg2.12098},
url = {https://ietresearch.onlinelibrary.wiley.com/doi/abs/10.1049/rpg2.12098},
eprint = {https://ietresearch.onlinelibrary.wiley.com/doi/pdf/10.1049/rpg2.12098},
year = {2021}
}

@article{OverflowMILP,
author = {Litlabø, Tormod and Aaslid, Per and Riise, Tarjei and Haugland, Dag},
year = {2023},
month = {08},
pages = {1-13},
title = {Modelling overflow using mixed integer programming in short-term hydropower scheduling},
journal = {Energy Systems},
doi = {10.1007/s12667-023-00602-2}
}

@article{NLPOpt,
author = {Arunkumar, R. and Jothiprakash, Vinayakam},
year = {2012},
month = {05},
pages = {111-120},
title = {Optimal Reservoir Operation for Hydropower Generation Using Non-Linear Programming Model},
volume = {93},
journal = {Journal of the Institution of Engineers (India): Serie A},
doi = {10.1007/s40030-012-0013-8}
}

@mastersthesis{AnderssonJohansson2025,
  author       = {Alfred Andersson and Holger Johansson},
  title        = {Hydropower optimization with detailed reservoir representation: Novel approaches for accurate modeling of reservoir levels in hydropower river systems},
  school       = {Chalmers University of Technology, Department of Space, Earth and Environment},
  year         = {2025},
  type         = {Master's thesis},
  address      = {Gothenburg, Sweden},
  url          = {https://www.chalmers.se/} 
}

@article{MultiMohan,
  author = {Mohan, S. and Raipure, Diwakar M.},
  title = {Multiobjective Analysis of Multireservoir System},
  journal = {Journal of Water Resources Planning and Management},
  volume = {118},
  number = {4},
  pages = {356--370},
  year = {1992},
  doi = {{10.1061/(ASCE)0733-9496(1992)118:4(356)}},
  url = {{https://ascelibrary.org/doi/abs/10.1061/(ASCE)0733-9496(1992)118:4(356)}},
  eprint = {{https://ascelibrary.org/doi/pdf/10.1061/(ASCE)0733-9496(1992)118:4(356)}}
}

@article{InflowInfluence,
author = {Barbosa, Alan and Celeste, Alcigeimes and Mendes, Ludmilson},
year = {2021},
month = {05},
pages = {},
title = {Influence of Inflow Nonstationarity on the Multipurpose Optimal Operation of Hydropower Plants Using Nonlinear Programming},
volume = {35},
journal = {Water Resources Management},
doi = {10.1007/s11269-021-02812-8}
}

@article{NIU2021107315,
title = {Multi-strategy gravitational search algorithm for constrained global optimization in coordinative operation of multiple hydropower reservoirs and solar photovoltaic power plants},
journal = {Applied Soft Computing},
volume = {107},
pages = {107315},
year = {2021},
issn = {1568-4946},
doi = {https://doi.org/10.1016/j.asoc.2021.107315},
url = {https://www.sciencedirect.com/science/article/pii/S1568494621002386},
author = {Wen-jing Niu and Zhong-kai Feng and Shuai Liu}
}

@article{WU2024121502,
title = {Short-Term Hydro-Wind-PV peak shaving scheduling using approximate hydropower output characters},
journal = {Renewable Energy},
volume = {236},
pages = {121502},
year = {2024},
issn = {0960-1481},
doi = {https://doi.org/10.1016/j.renene.2024.121502},
url = {https://www.sciencedirect.com/science/article/pii/S0960148124015702},
author = {Xinyu Wu and Jiaao Zhang and Xingchen Wei and Chuntian Cheng and Ruixiang Cheng}
}

@article{RO-Xu,
author = {Xu, Bin and Sun, Yu and Huang, Xin and Zhong, Ping-an and Zhu, Feilin and Zhang, Jian-yun and Wang, Xiaojun and Wang, Guoqing and Ma, Yufei and Qingwen, Lu and Wang, Han and Guo, Le},
year = {2022},
month = {04},
pages = {},
title = {Scenario‐Based Multiobjective Robust Optimization and Decision‐Making Framework for Optimal Operation of a Cascade Hydropower System Under Multiple Uncertainties},
volume = {58},
journal = {Water Resources Research},
doi = {10.1029/2021WR030965}
}

@ARTICLE{9887899,
  author={Wang, Chong and Ju, Ping and Wan, Can and Wu, Feng and Lei, Shunbo and Pan, Xueping and Lu, Tianguang},
  journal={IEEE Transactions on Sustainable Energy}, 
  title={Resilience-Based Coordinated Scheduling of Cascaded Hydro Power With Sequential Heavy Precipitation}, 
  year={2023},
  volume={14},
  number={2},
  pages={1299-1311},
  doi={10.1109/TSTE.2022.3205688}}

@article{AMINABADI2024282,
title = {Hybrid Genetic Algorithms and Heuristics for Nonlinear Short-Term Hydropower Optimization: A Comparative Analysis},
journal = {Procedia Computer Science},
volume = {246},
pages = {282-291},
year = {2024},
note = {28th International Conference on Knowledge Based and Intelligent information and Engineering Systems (KES 2024)},
issn = {1877-0509},
doi = {https://doi.org/10.1016/j.procs.2024.09.405},
url = {https://www.sciencedirect.com/science/article/pii/S1877050924024451},
author = {M. Jafari Aminabadi and S. Séguin and I. Fofana}
}

@INPROCEEDINGS{10353562,
  author={Wei, MingKui and Lu, Liang and Zhou, Hong and Shen, Li and Wang, Qing and Hu, Bangan},
  booktitle={2023 8th International Conference on Power and Renewable Energy (ICPRE)}, 
  title={Optimization Scheduling and Adequacy Evaluation of Cascade Hydropower Stations Considering the Extreme Weather}, 
  year={2023},
  volume={},
  number={},
  pages={809-814},
  doi={10.1109/ICPRE59655.2023.10353562}}

@article{DOGANIS20142051,
title = {Optimization of power production through coordinated use of hydroelectric and conventional power units},
journal = {Applied Mathematical Modelling},
volume = {38},
number = {7},
pages = {2051-2062},
year = {2014},
issn = {0307-904X},
doi = {https://doi.org/10.1016/j.apm.2013.10.025},
url = {https://www.sciencedirect.com/science/article/pii/S0307904X13006379},
author = {Philip Doganis and Haralambos Sarimveis}
}

@phdthesis{Forknall2014,
  author      = {Clayton Forknall},
  title       = {Optimal Management Strategies for a Cascade Reservoir System},
  school      = {University of Southern Queensland},
  institution = {Faculty of Health, Engineering and Sciences},
  year        = {2014},
  month       = {September},
  type        = {B.Sc. (Honours) Thesis},
  address     = {Toowoomba, Australia}
}

@article{YIN2022114582,
title = {Water-energy-ecosystem nexus modeling using multi-objective, non-linear programming in a regulated river: Exploring tradeoffs among environmental flows, cascaded small hydropower, and inter-basin water diversion projects},
journal = {Journal of Environmental Management},
volume = {308},
pages = {114582},
year = {2022},
issn = {0301-4797},
doi = {https://doi.org/10.1016/j.jenvman.2022.114582},
url = {https://www.sciencedirect.com/science/article/pii/S0301479722001554},
author = {Dongqin Yin and Xiang Li and Fang Wang and Yang Liu and Barry F.W. Croke and Anthony J. Jakeman}
}

@article{Wang2017OptOutputError,
  title   = {Study on optimization of the short-term operation of cascade hydropower stations by considering output error},
  author  = {Wang, Liping and Wang, Boquan and Zhang, Pu and Liu, Minghao and Li, Chuangang},
  journal = {Journal of Hydrology},
  volume  = {549},
  pages   = {326-339},
  year    = {2017},
  doi     = {10.1016/j.jhydrol.2017.03.074},
  url     = {https://doi.org/10.1016/j.jhydrol.2017.03.074},
  publisher = {Elsevier B.V.}
}

@article{Liu2025DynamicPS,
  title        = {Study on the peak shaving operation of cascade hydropower stations based on the plant-wide optimal curve},
  author       = {Liu, Fengshuo and Huang, Kui and Shi, Xuanyu and Zhao, Longqing and Yu, Yangxin and Ai, Xueshan and Fu, Xiang},
  journal      = {International Journal of Electrical Power and Energy Systems},
  volume       = {170},
  number       = {},
  pages        = {110920},
  year         = {2025},
  doi          = {10.1016/j.ijepes.2025.110920},
  url          = {https://doi.org/10.1016/j.ijepes.2025.110920},
  publisher    = {Elsevier},
  note         = {Open access under CC BY-NC-ND 4.0. Available online 30 July 2025}
}

@article{Zhu2025PFMODO,
  title        = {Optimizing multi-objective operation of cascade reservoirs to enhance hydropower generation and fish spawning demand synergies during the pre-flood drawdown period},
  author       = {Zhu, Di and Wang, Lin and Guo, Wei and Chen, Junhong and Li, Liping and Bu, Hui and Zuo, Jian and Chen, Hua},
  journal      = {Journal of Hydrology: Regional Studies},
  year         = {2025},
  volume       = {62},
  pages        = {102855},
  doi          = {10.1016/j.ejrh.2025.102855},
  url          = {https://doi.org/10.1016/j.ejrh.2025.102855},
  issn         = {2214-5818},
  publisher    = {Elsevier},
  note         = {Open access under CC BY 4.0; Available online 14 October 2025}
}

@article{NIU2018562,
title = {A parallel multi-objective particle swarm optimization for cascade hydropower reservoir operation in southwest China},
journal = {Applied Soft Computing},
volume = {70},
pages = {562-575},
year = {2018},
issn = {1568-4946},
doi = {https://doi.org/10.1016/j.asoc.2018.06.011},
url = {https://www.sciencedirect.com/science/article/pii/S156849461830334X},
author = {Wen-jing Niu and Zhong-kai Feng and Chun-tian Cheng and Xin-yu Wu}
}

@article{LIAO2023127685,
title = {Solution framework for short-term cascade hydropower system optimization operations based on the load decomposition strategy},
journal = {Energy},
volume = {277},
pages = {127685},
year = {2023},
issn = {0360-5442},
doi = {https://doi.org/10.1016/j.energy.2023.127685},
url = {https://www.sciencedirect.com/science/article/pii/S0360544223010794},
author = {Shengli Liao and Huan Liu and Benxi Liu and Tian Liu and Chonghao Li and Huaying Su}
}

@inproceedings{Liu2016_JointOptimization,
  author    = {Liu, Yuan and Zhou, Jianzhong and Chang, Chuyang and Lu, Peng and Wang, Chao and Tayyab, Muhammad},
  title     = {Short-Term Joint Optimization of Cascade Hydropower Stations on Daily Power Load Curve},
  booktitle = {2016 IEEE International Conference on Knowledge Engineering and Applications (ICKEA)},
  year = {2016},
  pages     = {236-240},
  publisher = {IEEE},
  isbn      = {978-1-5090-3471-0},
}

@article{xiao2023iboa,
  title = {An improved butterfly optimization algorithm and its application in cascade hydropower generation operation},
  author = {Xiao, Zhangling and Liang, Zhongmin and Wang, Jian and Li, Binquan and Hu, Yiming and Wang, Jun},
  journal = {Journal of Hydroinformatics},
  volume = {25},
  number = {3},
  pages = {1121--1138},
  year = {2023},
  doi = {10.2166/hydro.2023.026},
  url = {https://iwaponline.com/jh/article-pdf/25/3/1121/1228309/jh0251121.pdf},
  publisher = {IWA Publishing}
}

@article{Tayebiyan2016OptCleanEnergy,
  author  = {Tayebiyan, Aida and Mohammad, Thamer Ahmad},
  title   = {Optimization of cascade hydropower system operation by genetic algorithm to maximize clean energy output},
  journal = {Environmental Health Engineering and Management Journal},
  year    = {2016},
  volume  = {3},
  number  = {2},
  pages   = {99-106},
  doi     = {10.15171/EHEM.2016.07},
  note    = {ePublished: 18 June 2016}
}

@article{ManikkuwahandiHornberger2021,
  author    = {Thushara De Silva Manikkuwahandi and George M. Hornberger},
  title     = {Deriving Reservoir Cascade Operation Rules for Variable Streamflows by Optimizing Hydropower Generation and Irrigation Water Delivery},
  journal   = {Journal of Water Resources Planning and Management},
  year      = {2021},
  volume    = {147},
  number    = {7},
  pages     = {05021007},
  doi       = {10.1061/(ASCE)WR.1943-5452.0001372},
  issn      = {0733-9496},
  publisher = {American Society of Civil Engineers},
  url       = {https://doi.org/10.1061/(ASCE)WR.1943-5452.0001372}
}

@article{Chen2023CascadeHydropower,
  title        = {Cascade Hydropower System Operation Considering Ecological Flow Based on Different Multi-Objective Genetic Algorithms},
  author       = {Chen, Yubin and Wang, Manlin and Zhang, Yu and Lu, Yan and Xu, Bin and Yu, Lei},
  journal      = {Water Resources Management},
  year         = {2023},
  volume       = {37},
  number       = {10},
  pages        = {3093-3110},
  doi          = {10.1007/s11269-023-03491-3},
  url          = {https://doi.org/10.1007/s11269-023-03491-3},
  publisher    = {Springer Nature B.V.},
  note         = {Published online: 24 March 2023}
}

@article{Zhang2019Cascade,
  author       = {Hongxue Zhang and Jianxia Chang and Chao Gao and Hongshi Wu and Yimin Wang and Kaixuan Lei and Ruihao Long and Lianpeng Zhang},
  title        = {Cascade hydropower plants operation considering comprehensive ecological water demands},
  journal      = {Energy Conversion and Management},
  year         = {2019},
  volume       = {180},
  pages        = {119-133},
  doi          = {10.1016/j.enconman.2018.10.072},
  url          = {https://doi.org/10.1016/j.enconman.2018.10.072},
  note         = {Available online: 05 November 2018}
}

@article{Niu2021CSA, title = {Cooperation Search Algorithm for Power Generation Production Operation Optimization of Cascade Hydropower Reservoirs}, author = {Niu, Wen-jing and Feng, Zhong-kai and Li, Yu-rong and Liu, Shuai}, journal = {Water Resources Management}, year = {2021}, volume = {35}, number = {}, pages = {2465–2485}, publisher = {Springer Nature}, doi = {10.1007/s11269-021-02842-2}, url = {https://doi.org/10.1007/s11269-021-02842-2}, note = {Published online: 14 May 2021} }

@article{Sharifi2021FDBMSA, author = {Sharifi, Mohammad Reza and Akbarifard, Saeid and Qaderi, Kourosh and Madadi, Mohamad Reza}, title = {Developing MSA Algorithm by New Fitness-Distance-Balance Selection Method to Optimize Cascade Hydropower Reservoirs Operation}, journal = {Water Resources Management}, year = {2021}, volume = {35}, pages = {385–406}, doi = {10.1007/s11269-020-02745-8}, url = {https://doi.org/10.1007/s11269-020-02745-8}, publisher = {Springer Nature B.V.}, note = {Published online: 3 January 2021} }

@article{Li2025CascadeHydroIntervals,
  title        = {An Optimization Method for Day-Ahead Generation Interval of Cascade Hydropower Adapting to Multi-Source Coordinated Scheduling Requirements},
  author       = {Li, Shushan and Li, Chonghao and Wu, Huijun and Zhao, Zhipeng and Wang, Huan and Kang, Yongxi and Cheng, Chuntian and Li, Changhong},
  journal      = {Energies},
  year         = {2025},
  volume       = {18},
  pages        = {4901},
  doi          = {10.3390/en18184901},
  url          = {https://doi.org/10.3390/en18184901},
  publisher    = {MDPI},
  month        = {September}
}

@Article{ShangNested,
AUTHOR = {Shang, Ling and Li, Xiaofei and Shi, Haifeng and Kong, Feng and Wang, Ying and Shang, Yizi},
TITLE = {Long-, Medium-, and Short-Term Nested Optimized-Scheduling Model for Cascade Hydropower Plants: Development and Practical Application},
JOURNAL = {Water},
VOLUME = {14},
YEAR = {2022},
NUMBER = {10},
ARTICLE-NUMBER = {1586},
URL = {https://www.mdpi.com/2073-4441/14/10/1586},
ISSN = {2073-4441},
DOI = {10.3390/w14101586}
}

@Article{en18040964,
AUTHOR = {Wei, Daohong and Feng, Chunpeng and Liu, Dong},
TITLE = {Research on Economic Operation of Cascade Small Hydropower Stations Within Plants Based on Refined Efficiency Models},
JOURNAL = {Energies},
VOLUME = {18},
YEAR = {2025},
NUMBER = {4},
ARTICLE-NUMBER = {964},
URL = {https://www.mdpi.com/1996-1073/18/4/964},
ISSN = {1996-1073},
DOI = {10.3390/en18040964}
}

@article{JointOpt,
author = {Xiao Li  and Pan Liu  and Bo Ming  and Kangdi Huang  and Weifeng Xu  and Yan Wen },
title = {Joint Optimization of Forward Contract and Operating Rules for Cascade Hydropower Reservoirs},
journal = {Journal of Water Resources Planning and Management},
volume = {148},
number = {2},
pages = {04021099},
year = {2022},
doi = {10.1061/(ASCE)WR.1943-5452.0001510},

URL = {https://ascelibrary.org/doi/abs/10.1061/(ASCE)WR.1943-5452.0001510},
eprint = {https://ascelibrary.org/doi/pdf/10.1061/(ASCE)WR.1943-5452.0001510}
}

@Article{su152416916,
AUTHOR = {Li, Yang and Wu, Feng and Song, Xudong and Shi, Linjun and Lin, Keman and Hong, Feilong},
TITLE = {Data-Driven Chance-Constrained Schedule Optimization of Cascaded Hydropower and Photovoltaic Complementary Generation Systems for Shaving Peak Loads},
JOURNAL = {Sustainability},
VOLUME = {15},
YEAR = {2023},
NUMBER = {24},
ARTICLE-NUMBER = {16916},
URL = {https://www.mdpi.com/2071-1050/15/24/16916},
ISSN = {2071-1050},
DOI = {10.3390/su152416916}
}

@ARTICLE{627890,
  author={Xiaohong Guan and Ernan Ni and Renhou Li and Luh, P.B.},
  journal={IEEE Transactions on Power Systems}, 
  title={An optimization-based algorithm for scheduling hydrothermal power systems with cascaded reservoirs and discrete hydro constraints}, 
  year={1997},
  volume={12},
  number={4},
  pages={1775-1780},
  doi={10.1109/59.627890}}

@article{QIU2021124239,
title = {Risk analysis of water supply-hydropower generation-environment nexus in the cascade reservoir operation},
journal = {Journal of Cleaner Production},
volume = {283},
pages = {124239},
year = {2021},
issn = {0959-6526},
doi = {https://doi.org/10.1016/j.jclepro.2020.124239},
url = {https://www.sciencedirect.com/science/article/pii/S0959652620342840},
author = {Hongya Qiu and Lu Chen and Jianzhong Zhou and Zhongzheng He and Hansong Zhang}
}

@INPROCEEDINGS{11234414,
  author={Liu, Zhanzhi and Xu, Zhanxing and Li, Menglu and Wu, Haowei and Liu, Qiang and Li, Yang},
  booktitle={2025 9th International Conference on Power Energy Systems and Applications (ICoPESA)}, 
  title={Day-ahead Peak Shaving Dispatch of Photovoltaic-Hybrid Pumped Storage Hydropower Systems}, 
  year={2025},
  volume={},
  number={},
  pages={572-576},
  doi={10.1109/ICoPESA65876.2025.11234414}}

@INPROCEEDINGS{10672978,
  author={Wu, Zhengqiang and Fan, Jiesheng and Zhang, Quan’e and Han, Huan and Song, Juwen and Jin, Zhongxia},
  booktitle={2024 7th International Conference on Electronics Technology (ICET)}, 
  title={Optimization Scheduling for Cascaded Hydropower Groups Adapted to Tiered Carbon Trading}, 
  year={2024},
  volume={},
  number={},
  pages={501-506},
  doi={10.1109/ICET61945.2024.10672978}}

@INPROCEEDINGS{7853589,
  author={Sharma, Prakash Chand and Abhyankar, A.R.},
  booktitle={2016 IEEE 1st International Conference on Power Electronics, Intelligent Control and Energy Systems (ICPEICES)}, 
  title={Multi-objective short-term tandem hydro scheduling: MINLP approach}, 
  year={2016},
  volume={},
  number={},
  pages={1-6},
  doi={10.1109/ICPEICES.2016.7853589}}

@INPROCEEDINGS{11009646,
  author={Li, Chengyuan and Ma, Luyu and Yu, Xiaowei},
  booktitle={2025 2nd International Conference on Smart Grid and Artificial Intelligence (SGAI)}, 
  title={Long-Term Reversible Hydropower System Operation Optimization Based on Fast Unit Commitment}, 
  year={2025},
  volume={},
  number={},
  pages={791-794},
  doi={10.1109/SGAI64825.2025.11009646}}

@INPROCEEDINGS{4839958,
  author={Wenping Chang and Xianjue Luo and Hai Yu},
  booktitle={2009 IEEE/PES Power Systems Conference and Exposition}, 
  title={A fuzzy adaptive particle swarm optimization for Long-Term Optimal Scheduling of Cascaded hydropower station}, 
  year={2009},
  volume={},
  number={},
  pages={1-5},
  doi={10.1109/PSCE.2009.4839958}}

@INPROCEEDINGS{5287657,
  author={Lihua, Chen and Yadong, Mei and Na, Yang},
  booktitle={2009 Second International Conference on Intelligent Computation Technology and Automation}, 
  title={Parallel Particle Swarm Optimization Algorithm and Its Application in the Optimal Operation of Cascade Reservoirs in Yalong River}, 
  year={2009},
  volume={1},
  number={},
  pages={279-282},
  doi={10.1109/ICICTA.2009.75}}

@INPROCEEDINGS{4601846,
  author={Moosavian, S. Ali A. and Ghafari, A. and Salimi, A. and Abdi, N.},
  booktitle={2008 IEEE/ASME International Conference on Advanced Intelligent Mechatronics}, 
  title={Non-linear multiobjective optimization for control of hydropower plants network}, 
  year={2008},
  volume={},
  number={},
  pages={1278-1283},
  doi={10.1109/AIM.2008.4601846}}

@INPROCEEDINGS{964906,
  author={Garcia-Gonzalez, J. and Castro, G.A.},
  booktitle={2001 IEEE Porto Power Tech Proceedings (Cat. No.01EX502)}, 
  title={Short-term hydro scheduling with cascaded and head-dependent reservoirs based on mixed-integer linear programming}, 
  year={2001},
  volume={3},
  number={},
  pages={6 pp. vol.3-},
  doi={10.1109/PTC.2001.964906}}

@INPROCEEDINGS{6009930,
  author={Jiangtao Jia and Xiaohong Guan},
  booktitle={2011 2nd International Conference on Artificial Intelligence, Management Science and Electronic Commerce (AIMSEC)}, 
  title={MILP formulation for short-term scheduling of cascaded reservoirs with head effects}, 
  year={2011},
  volume={},
  number={},
  pages={4061-4064},
  doi={10.1109/AIMSEC.2011.6009930}}

@INPROCEEDINGS{7443262,
  author={Sutradhar, Suman and Choudhury, N B Dev and Sinha, N},
  booktitle={2015 Annual IEEE India Conference (INDICON)}, 
  title={Mixed integer non-linear programming for hydrothermal scheduling problem}, 
  year={2015},
  volume={},
  number={},
  pages={1-6},
  doi={10.1109/INDICON.2015.7443262}}

@INPROCEEDINGS{7804914,
  author={Bao, Gang and Wen, Siyu and Yu, Wenju},
  booktitle={2016 31st Youth Academic Annual Conference of Chinese Association of Automation (YAC)}, 
  title={Neural network approach for modelling and solving the unit commitment problem of cascaded hydropower stations}, 
  year={2016},
  volume={},
  number={},
  pages={334-338},
  doi={10.1109/YAC.2016.7804914}}

@INPROCEEDINGS{5953023,
  author={Eusébio, Eduardo and Camus, Cristina and Mendes, Victor},
  booktitle={2011 8th International Conference on the European Energy Market (EEM)}, 
  title={Short-term value for the water stored in head-sensitivity power system reservoirs}, 
  year={2011},
  volume={},
  number={},
  pages={275-280},
  doi={10.1109/EEM.2011.5953023}}

@INPROCEEDINGS{7853419,
  author={Sutradhar, Suman and Choudhury, N.B. Dev and Sinha, N.},
  booktitle={2016 IEEE 1st International Conference on Power Electronics, Intelligent Control and Energy Systems (ICPEICES)}, 
  title={MINLP for Hydro-Thermal Unit Commitment problem using BONMIN solver}, 
  year={2016},
  volume={},
  number={},
  pages={1-6},
  doi={10.1109/ICPEICES.2016.7853419}}

@INPROCEEDINGS{6019203,
  author={Pousinho, H. M. I. and Mendes, V. M. F. and Catalão, J. P. S.},
  booktitle={2011 IEEE Trondheim PowerTech}, 
  title={Profit-based head-dependent short-term hydro scheduling considering risk constraints}, 
  year={2011},
  volume={},
  number={},
  pages={1-6},
  doi={10.1109/PTC.2011.6019203}}

@INPROCEEDINGS{6072098,
  author={Kladnik, B. and Artač, G. and Kozan, B. and Gubina, A. F. and Nagode, K. and Dusak, M.},
  booktitle={IEEE Africon '11}, 
  title={Scheduling the Slovenian cascaded hydro system on the river Sava}, 
  year={2011},
  volume={},
  number={},
  pages={1-5},
  doi={10.1109/AFRCON.2011.6072098}}

@INPROCEEDINGS{7483057,
  author={Heidarizadeh, Mohammad and Shivaie, Mojtaba and Ahmadian, Mohammad and Ameli, Mohammad T.},
  booktitle={2016 6th Conference on Thermal Power Plants (CTPP)}, 
  title={A risk-based optimal self-scheduling of cascaded hydro power plants in joint energy and reserve electricity markets}, 
  year={2016},
  volume={},
  number={},
  pages={76-82},
  doi={10.1109/CTPP.2016.7483057}}

@INPROCEEDINGS{7838997,
  author={Sutradhar, Suman and Dev Choudhury, N. B. and Sinha, N},
  booktitle={2016 IEEE Annual India Conference (INDICON)}, 
  title={MINLP for hydrothermal scheduling problem considering transmission loss}, 
  year={2016},
  volume={},
  number={},
  pages={1-6},
  doi={10.1109/INDICON.2016.7838997}}

@article{LIU202311,
title = {Short-term operation of cascade hydropower system sharing flexibility via high voltage direct current lines for multiple grids peak shaving},
journal = {Renewable Energy},
volume = {213},
pages = {11-29},
year = {2023},
issn = {0960-1481},
doi = {https://doi.org/10.1016/j.renene.2023.05.095},
url = {https://www.sciencedirect.com/science/article/pii/S0960148123007279},
author = {Benxi Liu and Tengyuan Liu and Shengli Liao and Haidong Wang and Xiaoyu Jin}
}

@article{JIN2024132283,
title = {Assessing hydropower capability for accommodating variable renewable energy considering peak shaving of multiple power grids},
journal = {Energy},
volume = {305},
pages = {132283},
year = {2024},
issn = {0360-5442},
doi = {https://doi.org/10.1016/j.energy.2024.132283},
url = {https://www.sciencedirect.com/science/article/pii/S0360544224020577},
author = {Xiaoyu Jin and Benxi Liu and Shengli Liao and Chuntian Cheng and Yi Zhang and Zebin Jia}
}

@article{ZHANG2019883,
title = {Coordinated optimal operation of hydro–wind–solar integrated systems},
journal = {Applied Energy},
volume = {242},
pages = {883-896},
year = {2019},
issn = {0306-2619},
doi = {https://doi.org/10.1016/j.apenergy.2019.03.064},
url = {https://www.sciencedirect.com/science/article/pii/S0306261919304738},
author = {Hongxuan Zhang and Zongxiang Lu and Wei Hu and Yiting Wang and Ling Dong and Jietan Zhang}
}

@article{MA2023138602,
title = {Multi-objective solution and decision-making framework for coordinating the short-term hydropeaking-navigation-production conflict of cascade hydropower reservoirs},
journal = {Journal of Cleaner Production},
volume = {422},
pages = {138602},
year = {2023},
issn = {0959-6526},
doi = {https://doi.org/10.1016/j.jclepro.2023.138602},
url = {https://www.sciencedirect.com/science/article/pii/S0959652623027609},
author = {Xiangyu Ma and Shengli Liao and Benxi Liu and Hongye Zhao and Chuntian Cheng and Huaying Su}
}

@INPROCEEDINGS{10853341,
  author={Tao, Haiguo and Guo, Jiang and Liu, Wei and Chen, Ci and Li, Lianjie and Chen, Yanli and Zong, Wei},
  booktitle={2nd International Conference on Power, Communication, Computing and Networking Technologies (PCCNT 2024)}, 
  title={Research on medium and long term optimal operation of cascade hydropower stations on Tibet's Lhasa River-Yarlung Zangbo River}, 
  year={2024},
  volume={2024},
  number={},
  pages={263-273},
  doi={10.1049/icp.2024.4116}}

@article{LIU2023127298,
title = {Effect of the quality of streamflow forecasts on the operation of cascade hydropower stations using stochastic optimization models},
journal = {Energy},
volume = {273},
pages = {127298},
year = {2023},
issn = {0360-5442},
doi = {https://doi.org/10.1016/j.energy.2023.127298},
url = {https://www.sciencedirect.com/science/article/pii/S0360544223006928},
author = {Yuan Liu and Changming Ji and Yi Wang and Yanke Zhang and Zhiqiang Jiang and Qiumei Ma and Xiaoning Hou}
}

@Article{su151310002,
AUTHOR = {Coban, Hasan Huseyin},
TITLE = {Hydropower Planning in Combination with Batteries and Solar Energy},
JOURNAL = {Sustainability},
VOLUME = {15},
YEAR = {2023},
NUMBER = {13},
ARTICLE-NUMBER = {10002},
URL = {https://www.mdpi.com/2071-1050/15/13/10002},
ISSN = {2071-1050},
DOI = {10.3390/su151310002}
}

@article{LU202425,
title = {Stochastic programming based coordinated expansion planning of generation, transmission, demand side resources, and energy storage considering the DC transmission system},
journal = {Global Energy Interconnection},
volume = {7},
number = {1},
pages = {25-37},
year = {2024},
issn = {2096-5117},
doi = {https://doi.org/10.1016/j.gloei.2024.01.003},
url = {https://www.sciencedirect.com/science/article/pii/S2096511724000033},
author = {Liang Lu and Mingkui Wei and Yuxuan Tao and Qing Wang and Yuxiao Yang and Chuan He and Haonan Zhang}
}

@article{LU2022435,
title = {Day-ahead optimal dispatching of multi-source power system},
journal = {Renewable Energy},
volume = {183},
pages = {435-446},
year = {2022},
issn = {0960-1481},
doi = {https://doi.org/10.1016/j.renene.2021.10.093},
url = {https://www.sciencedirect.com/science/article/pii/S0960148121015597},
author = {Mengke Lu and Jun Guan and Huahua Wu and Huizhe Chen and Wei Gu and Ye Wu and ChengXiang Ling and Linqiang Zhang}
}

@article{ZHANG2023108833,
title = {A mathematical programming–based heuristic for coordinated hydrothermal generator maintenance scheduling and long-term unit commitment},
journal = {International Journal of Electrical Power \& Energy Systems},
volume = {147},
pages = {108833},
year = {2023},
issn = {0142-0615},
doi = {https://doi.org/10.1016/j.ijepes.2022.108833},
url = {https://www.sciencedirect.com/science/article/pii/S0142061522008298},
author = {Zihan Zhang and Mingbo Liu and Min Xie and Ping Dong}
}

@article{Knezevic2016234245,
author = {Knežević, Goran and Baus, Zoran and Nikolovski, Srete},
doi = {10.1515/jee-2016-0035},
url = {https://doi.org/10.1515/jee-2016-0035},
title = {Short–Term Planning of Hybrid Power System},
journal = {Journal of Electrical Engineering},
number = {4},
volume = {67},
year = {2016},
pages = {234--245}
}

@article{Xu2017multiobjective,
  author    = {Bin Xu and Ping-an Zhong and Yenan Wu and Fangming Fu and Yuting Chen and Yunfa Zhao},
  title     = {A Multiobjective Stochastic Programming Model for Hydropower Hedging Operations under Inexact Information},
  journal   = {Water Resources Management},
  year      = {2017},
  volume    = {31},
  number    = {15-16},
  pages     = {4649--4667},
  doi       = {10.1007/s11269-017-1771-x},
  url       = {https://doi.org/10.1007/s11269-017-1771-x},
  publisher = {Springer Science+Business Media B.V.}
}

@article{KANG2026124031,
title = {Integrating water delay time into short-term hydropower scheduling with spinning reserve capacity allocation and execution},
journal = {Renewable Energy},
volume = {256},
pages = {124031},
year = {2026},
issn = {0960-1481},
doi = {https://doi.org/10.1016/j.renene.2025.124031},
url = {https://www.sciencedirect.com/science/article/pii/S0960148125016957},
author = {Yongxi Kang and Zhipeng Zhao and Chuntian Cheng and Xiangyu Wu and Xiaoyu Jin and Huaying Su}
}

@article{LIU2025116775,
title = {Short-term complementary scheduling of cascade energy storage systems for wind and solar regulation},
journal = {Journal of Energy Storage},
volume = {124},
pages = {116775},
year = {2025},
issn = {2352-152X},
doi = {https://doi.org/10.1016/j.est.2025.116775},
url = {https://www.sciencedirect.com/science/article/pii/S2352152X25014884},
author = {Yuanyuan Liu and Hao Zhang and Pengcheng Guo and Shuai Wu}
}

@article{WANG2024130258,
title = {Multi-objective day-ahead scheduling of cascade hydropower-photovoltaic complementary system with pumping installation},
journal = {Energy},
volume = {290},
pages = {130258},
year = {2024},
issn = {0360-5442},
doi = {https://doi.org/10.1016/j.energy.2024.130258},
url = {https://www.sciencedirect.com/science/article/pii/S036054422400029X},
author = {Zizhao Wang and Yang Li and Feng Wu and Jiawei Wu and Linjun Shi and Keman Lin}
}

@article{HUANG2025132756,
title = {Framework for short-term hydropower cascade–station–unit integrated multi-objective scheduling: Considering unit safety and economic efficiency},
journal = {Journal of Hydrology},
volume = {653},
pages = {132756},
year = {2025},
issn = {0022-1694},
doi = {https://doi.org/10.1016/j.jhydrol.2025.132756},
url = {https://www.sciencedirect.com/science/article/pii/S0022169425000940},
author = {Jingwei Huang and Hui Qin and Xu Yang and Keyan Shen and Huaming Yao and Xinyu Chang and Gaoge Li and Yuan Gao}
}

@Article{en18143745,
AUTHOR = {Su, Huaying and Li, Yupeng and Zhang, Yan and Wang, Yujian and Li, Gang and Cheng, Chuntian},
TITLE = {A Mid-Term Scheduling Method for Cascade Hydropower Stations to Safeguard Against Continuous Extreme New Energy Fluctuations},
JOURNAL = {Energies},
VOLUME = {18},
YEAR = {2025},
NUMBER = {14},
ARTICLE-NUMBER = {3745},
URL = {https://www.mdpi.com/1996-1073/18/14/3745},
ISSN = {1996-1073},
DOI = {10.3390/en18143745}
}

@Article{w17101441,
AUTHOR = {Lu, Jia and Fang, Zhou and Zhang, Zheng and Liu, Yaxin and Xu, Yang and Wang, Tao and Yang, Yuqi},
TITLE = {Progressive Linear Programming Optimality Method Based on Decomposing Nonlinear Functions for Short-Term Cascade Hydropower Scheduling},
JOURNAL = {Water},
VOLUME = {17},
YEAR = {2025},
NUMBER = {10},
ARTICLE-NUMBER = {1441},
URL = {https://www.mdpi.com/2073-4441/17/10/1441},
ISSN = {2073-4441},
DOI = {10.3390/w17101441}
}

@article{LI2023103534,
title = {Risk-averse energy management of hydro/thermal/pumped storage complementarily operating with wind/solar: Balancing risk, cost and carbon emission},
journal = {Sustainable Energy Technologies and Assessments},
volume = {60},
pages = {103534},
year = {2023},
issn = {2213-1388},
doi = {https://doi.org/10.1016/j.seta.2023.103534},
url = {https://www.sciencedirect.com/science/article/pii/S2213138823005271},
author = {Xudong Li and Weijia Yang and Zhigao Zhao and Ran Wang and Xiuxing Yin and Pan Liu}
}

@Article{en17112734,
AUTHOR = {Liu, Yuanyuan and Zhang, Hao and Guo, Pengcheng and Li, Chenxi and Wu, Shuai},
TITLE = {Optimal Scheduling of a Cascade Hydropower Energy Storage System for Solar and Wind Energy Accommodation},
JOURNAL = {Energies},
VOLUME = {17},
YEAR = {2024},
NUMBER = {11},
ARTICLE-NUMBER = {2734},
URL = {https://www.mdpi.com/1996-1073/17/11/2734},
ISSN = {1996-1073},
DOI = {10.3390/en17112734}
}

@article{SU2025124360,
title = {Optimal scheduling of a cascade hydro-thermal-wind power system integrating data centers and considering the spatiotemporal asynchronous transfer of energy resources},
journal = {Applied Energy},
volume = {377},
pages = {124360},
year = {2025},
issn = {0306-2619},
doi = {https://doi.org/10.1016/j.apenergy.2024.124360},
url = {https://www.sciencedirect.com/science/article/pii/S0306261924017434},
author = {Chengguo Su and Lingshuang Wang and Quan Sui and Huijun Wu}
}

@article{Shen24,
author = {Shen, Li and Wang, Qing and Wan, Yizhi and Xu, Xiao and Liu, Youbo},
year = {2025},
month = {01},
pages = {},
title = {Multi-timescale scheduling optimization of cascade hydro-solar complementary power stations considering spatio-temporal correlation},
volume = {80},
journal = {Science and Technology for Energy Transition},
doi = {10.2516/stet/2024104}
}

@article{LU2021126388,
title = {Stochastic programming for floodwater utilization of a complex multi-reservoir system considering risk constraints},
journal = {Journal of Hydrology},
volume = {599},
pages = {126388},
year = {2021},
issn = {0022-1694},
doi = {https://doi.org/10.1016/j.jhydrol.2021.126388},
url = {https://www.sciencedirect.com/science/article/pii/S0022169421004352},
author = {Qingwen Lu and Ping-an Zhong and Bin Xu and Feilin Zhu and Xin Huang and Han Wang and Yufei Ma}
}

@article{Metamodeling,
author = {Piguet, Antoine and Benefice, Astrig and Bontron, Guillaume and Helbert, Celine and Vial, Grégory},
year = {2023},
month = {10},
pages = {},
title = {Metamodeling the optimal total revenues of the short-term optimization of a hydropower cascade under uncertainty},
volume = {78},
journal = {Science and Technology for Energy Transition},
doi = {10.2516/stet/2023026}
}

@article{MaxConsume,
author = {Hu, Shuzhe and Miao, Jinniu and Wu, Jingyang and Zhao, Liqian and Wang, Yue and Meng, Fanyan and Wei, Chao and Zhang, Xiaoqin and Zhu, Benrui},
year = {2025},
month = {01},
pages = {328},
title = {Short-Term Optimal Scheduling of a Cascade Hydro-Photovoltaic System for Maximizing the Expectation of Consumable Electricity},
volume = {13},
journal = {Processes},
doi = {10.3390/pr13020328}
}

@article{CHENG2022123908,
title = {A hierarchical model in short-term hydro scheduling with unit commitment and head-dependency},
journal = {Energy},
volume = {251},
pages = {123908},
year = {2022},
issn = {0360-5442},
doi = {https://doi.org/10.1016/j.energy.2022.123908},
url = {https://www.sciencedirect.com/science/article/pii/S0360544222008118},
author = {Xianliang Cheng and Suzhen Feng and Hao Zheng and Jinwen Wang and Shuangquan Liu}
}

@article{FENG2020119035,
title = {An effective three-stage hybrid optimization method for source-network-load power generation of cascade hydropower reservoirs serving multiple interconnected power grids},
journal = {Journal of Cleaner Production},
volume = {246},
pages = {119035},
year = {2020},
issn = {0959-6526},
doi = {https://doi.org/10.1016/j.jclepro.2019.119035},
url = {https://www.sciencedirect.com/science/article/pii/S0959652619339058},
author = {Zhong-kai Feng and Wen-jing Niu and Xiong Cheng and Jia-yang Wang and Sen Wang and Zhen-guo Song}
}

@ARTICLE{EFFLP,
       author = {{Kang}, Chuanxiong and {Chen}, Cheng and {Wang}, Jinwen},
        title = "{An Efficient Linearization Method for Long-Term Operation of Cascaded Hydropower Reservoirs}",
      journal = {Water Resources Management},
         year = 2018,
        month = aug,
       volume = {32},
       number = {10},
        pages = {3391-3404},
          doi = {10.1007/s11269-018-1997-2},
       adsurl = {https://ui.adsabs.harvard.edu/abs/2018WatRM..32.3391K}
}

@ARTICLE{8295134,
  author={Apostolopoulou, Dimitra and De Grève, Zacharie and McCulloch, Malcolm},
  journal={IEEE Transactions on Power Systems}, 
  title={Robust Optimization for Hydroelectric System Operation Under Uncertainty}, 
  year={2018},
  volume={33},
  number={3},
  pages={3337-3348},
  doi={10.1109/TPWRS.2018.2807794}}

@article{DELADURANTAYE2009499,
title = {Optimizing profits from hydroelectricity production},
journal = {Computers \& Operations Research},
volume = {36},
number = {2},
pages = {499-529},
year = {2009},
note = {Scheduling for Modern Manufacturing, Logistics, and Supply Chains},
issn = {0305-0548},
doi = {https://doi.org/10.1016/j.cor.2007.10.012},
url = {https://www.sciencedirect.com/science/article/pii/S0305054807002079},
author = {Daniel {De Ladurantaye} and Michel Gendreau and Jean-Yves Potvin}
}

@article{LargeScaleOpt,
  author = {Barros, Mario T. L. and Tsai, Frank T-C. and Yang, Shu-li and Lopes, Joao E. G. and Yeh, William W-G.},
  title = {Optimization of Large-Scale Hydropower System Operations},
  journal = {Journal of Water Resources Planning and Management},
  volume = {129},
  number = {3},
  pages = {178--188},
  year = {2003},
  doi = {10.1061/(ASCE)0733-9496(2003)129:3(178)},
  url = {https://ascelibrary.org/doi/abs/10.1061/(ASCE)0733-9496(2003)129:3(178)},
  eprint = {10.1061/(ASCE)0733-9496(2003)129:3(178)}
}

@article{PLOUSSARD202450,
title = {Piecewise linear approximation with minimum number of linear segments and minimum error: A fast approach to tighten and warm start the hierarchical mixed integer formulation},
journal = {European Journal of Operational Research},
volume = {315},
number = {1},
pages = {50-62},
year = {2024},
issn = {0377-2217},
doi = {https://doi.org/10.1016/j.ejor.2023.11.017},
url = {https://www.sciencedirect.com/science/article/pii/S0377221723008585},
author = {Quentin Ploussard}
}

@article{ploussard2025tightening,
  author  = {Ploussard, Quentin and Li, Xiang and Pavi{\v{c}}evi{\'c}, Matija},
  title   = {Tightening the Difference-of-Convex Formulation for the Piecewise Linear Approximation in General Dimensions},
  journal = {INFORMS Journal on Optimization},
  year    = {2025},
  doi     = {10.1287/ijoo.2025.0074},
  url     = {https://pubsonline.informs.org/doi/10.1287/ijoo.2025.0074}
}

@article{vielma2010nonseparable,
  author  = {Vielma, Juan Pablo and Ahmed, Shabbir and Nemhauser, George L.},
  title   = {Mixed-Integer Models for Nonseparable Piecewise Linear Optimization: Unifying Framework and Extensions},
  journal = {Operations Research},
  volume  = {58},
  number  = {2},
  pages   = {303--315},
  year    = {2010},
  doi     = {10.1287/opre.1090.0721},
  url     = {https://pubsonline.informs.org/doi/10.1287/opre.1090.0721}
}

@article{GRIMSTAD2019106580,
title = {ReLU networks as surrogate models in mixed-integer linear programs},
journal = {Computers \& Chemical Engineering},
volume = {131},
pages = {106580},
year = {2019},
issn = {0098-1354},
doi = {https://doi.org/10.1016/j.compchemeng.2019.106580},
url = {https://www.sciencedirect.com/science/article/pii/S0098135419307203},
author = {Bjarne Grimstad and Henrik Andersson}
}

@article{huang_relu_2020,
	title = {{ReLU} {Networks} {Are} {Universal} {Approximators} via {Piecewise} {Linear} or {Constant} {Functions}},
	volume = {32},
	issn = {0899-7667},
	url ={https://doi.org/10.1162/neco_a_01316},
	doi = {10.1162/neco_a_01316},
	number = {11},
	urldate = {2024-12-29},
	journal = {Neural Computation},
	author = {Huang, Changcun},
	month = nov,
	year = {2020},
	pages = {2249--2278},
}

@article{THIRUNAVUKKARASU2023113192,
title = {A comprehensive review on optimization of hybrid renewable energy systems using various optimization techniques},
journal = {Renewable and Sustainable Energy Reviews},
volume = {176},
pages = {113192},
year = {2023},
issn = {1364-0321},
doi = {https://doi.org/10.1016/j.rser.2023.113192},
url = {https://www.sciencedirect.com/science/article/pii/S1364032123000485},
author = {M. Thirunavukkarasu and Yashwant Sawle and Himadri Lala}
}

@article{taktak2017overview,
  author  = {Taktak, Raouia and D'Ambrosio, Claudia},
  title   = {An overview on mathematical programming approaches for the deterministic unit commitment problem in hydro valleys},
  journal = {Energy Systems},
  year    = {2017},
  volume  = {8},
  number  = {1},
  pages   = {57--79},
  doi     = {10.1007/s12667-015-0189-x},
  url     = {https://link.springer.com/article/10.1007/s12667-015-0189-x}
}

@article{labadie2004optimal,
  author  = {Labadie, John W.},
  title   = {Optimal Operation of Multireservoir Systems: State-of-the-Art Review},
  journal = {Journal of Water Resources Planning and Management},
  volume  = {130},
  number  = {2},
  pages   = {93--111},
  year    = {2004},
  doi     = {10.1061/(ASCE)0733-9496(2004)130:2(93)},
  url     = {https://ascelibrary.org/doi/10.1061/(ASCE)0733-9496(2004)130:2(93)}
}

@Article{w11071392,
AUTHOR = {Parvez, Iram and Shen, Jianjian and Khan, Mehran and Cheng, Chuntian},
TITLE = {Modeling and Solution Techniques Used for Hydro Generation Scheduling},
JOURNAL = {Water},
VOLUME = {11},
YEAR = {2019},
NUMBER = {7},
ARTICLE-NUMBER = {1392},
URL = {https://www.mdpi.com/2073-4441/11/7/1392},
ISSN = {2073-4441},
DOI = {10.3390/w11071392}
}

@article{lai2022review,
  author  = {Lai, Vivien and Huang, Yuk Feng and Koo, Chai Hoon and Ahmed, Ali Najah and El-Shafie, Ahmed},
  title   = {A Review of Reservoir Operation Optimisations: From Traditional Models to Metaheuristic Algorithms},
  journal = {Archives of Computational Methods in Engineering},
  volume  = {29},
  number  = {5},
  pages   = {3435--3457},
  year    = {2022},
  doi     = {10.1007/s11831-021-09701-8},
  url     = {https://link.springer.com/article/10.1007/s11831-021-09701-8}
}

@techreport{osti_1726280,
  author       = {Voisin, Nathalie and Bain, Dominique and Macknick, Jordan and O'Neil, Rebecca S.},
  title        = {Improving Hydropower Representation in Power System Models (Summary of Technical Workshop)},
  institution  = {Pacific Northwest National Laboratory (PNNL), Richland, WA (United States); National Laboratory of the Rockies (NLR), Golden, CO (United States)},
  doi          = {10.2172/1726280},
  url          = {https://www.osti.gov/biblio/1726280},
  place        = {United States},
  year         = {2020},
  month        = {11}}

@article{Huangfu2018HiGHS,
  author  = {Huangfu, Qi and Hall, Julian A. J.},
  title   = {Parallelizing the dual revised simplex method},
  journal = {Mathematical Programming Computation},
  year    = {2018},
  volume  = {10},
  number  = {1},
  pages   = {119--142},
  doi     = {10.1007/s12532-017-0130-5}
}

@article{Achterberg2009SCIP,
  author  = {Achterberg, Tobias},
  title   = {SCIP: Solving Constraint Integer Programs},
  journal = {Mathematical Programming Computation},
  year    = {2009},
  volume  = {1},
  number  = {1},
  pages   = {1--41},
  doi     = {10.1007/s12532-008-0001-1}
}

@book{williams_model_2013,
	title        = {Model building in mathematical programming},
	author       = {Williams, H. P.},
	location     = {Hoboken, N.J},
	publisher    = {Wiley},
	isbn         = {978-1-118-44333-0},
	edition      = {5th ed},
	pagetotal    = 411,
	date         = 2013,
    year = {2013}
}

@article{KAZDA2021107310,
title = {Nonconvex multivariate piecewise-linear fitting using the difference-of-convex representation},
journal = {Computers \& Chemical Engineering},
volume = {150},
pages = {107310},
year = {2021},
issn = {0098-1354},
doi = {https://doi.org/10.1016/j.compchemeng.2021.107310},
url = {https://www.sciencedirect.com/science/article/pii/S0098135421000880},
author = {Kody Kazda and Xiang Li}
}

@book{davis2006direct,
  title     = {Direct Methods for Sparse Linear Systems},
  author    = {Davis, Timothy A.},
  year      = {2006},
  publisher = {Society for Industrial and Applied Mathematics},
  address   = {Philadelphia, PA},
  doi       = {10.1137/1.9780898718881}
}

@techreport{osti_1922507,
  author = {Kincic, Slaven and Samaan, Nader and Datta, Sohom and Somani, Abhishek and Yuan, Haoyu and Tan, Jin and Bhattarai, Rojan and Mosier, Thomas M.},
  title = {Hydropower Modeling Gaps in Planning and Operational Studies},
  institution = {Pacific Northwest National Laboratory (PNNL), Richland, WA (United States); National Laboratory of the Rockies (NLR), Golden, CO (United States); Idaho National Laboratory (INL), Idaho Falls, ID (United States)},
  year = {2022},
  month = {11},
  doi = {10.2172/1922507},
  url = {https://www.osti.gov/biblio/1922507},
  place = {United States}
}

@book{junger2010fifty,
  editor    = {J{\"u}nger, Michael and Liebling, Thomas M. and Naddef, Denis and Nemhauser, George L. and Pulleyblank, William R. and Reinelt, Gerhard and Rinaldi, Giovanni and Wolsey, Laurence A.},
  title     = {50 Years of Integer Programming 1958--2008: From the Early Years to the State-of-the-Art},
  publisher = {Springer},
  address   = {Berlin, Heidelberg},
  year      = {2010},
  doi       = {10.1007/978-3-540-68279-0},
  isbn      = {978-3-540-68279-0}
}

\newpage

\begin{quote}
    \footnotesize
    The submitted manuscript has been created by UChicago Argonne, LLC, Operator of Argonne National Laboratory (“Argonne”). Argonne, a U.S. Department of Energy Office of Science laboratory, is operated under Contract No. DE-AC02-06CH11357. The U.S. Government retains for itself, and others acting on its behalf, a paid-up nonexclusive, irrevocable worldwide license in said article to reproduce, prepare derivative works, distribute copies to the public, and perform publicly and display publicly, by or on behalf of the Government.  The Department of Energy will provide public access to these results of federally sponsored research in accordance with the DOE Public Access Plan. http://energy.gov/downloads/doe-public-access-plan 
\end{quote}

\end{document}